\documentclass[12pt]{elsart}
\usepackage[T1]{fontenc}
\usepackage{amsmath, amssymb, amsfonts}
\usepackage{dsfont}
\usepackage{graphicx}
\usepackage{slashbox}
\usepackage{color}

\def\E{\mathbb{E}}
\def\pen{{\mathrm{pen}}}
\def\1{\mathds{1}}
\newcommand{\transpose}[1]{{\vphantom{#1}}^{\mathit t}{#1}}

\begin{document}

\begin{frontmatter}

\title{Adaptive estimation of the transition density of a Markov Chain} 

\author{Claire Lacour}


\address{Laboratoire MAP5, Universit\'e Paris 5, 
45, rue des Saints-P\`eres, 75270 Paris Cedex 06, France\\
lacour@math-info.univ-paris5.fr}

\begin{abstract}
In this paper a new estimator for the transition density $\pi$ of an homogeneous Markov chain is 
considered. We introduce an original contrast derived from regression framework and we use 
a model selection method to estimate $\pi$ under mild conditions. The resulting
estimate is  adaptive with an optimal rate of convergence over a large range
of anisotropic Besov spaces $B_{2,\infty}^{(\alpha_1,\alpha_2)}$. Some applications and simulations 
are also presented. 

\vspace{0.5cm}
\hspace{-0.5cm}\textbf{ R\'esum\'e } 

\hspace{-0.35cm}Dans cet article, on consid\`ere un nouvel estimateur de la densit\'e de transition 
$\pi$ d'une cha\^ine de Markov homog\`ene. Pour cela, on introduit un contraste original issu de la 
th\'eorie de la r\'egression et on utilise 
une m\'ethode de s\'election de mod\`eles pour estimer $\pi$ sous des conditions peu
restrictives. L'estimateur obtenu est adaptatif et la vitesse de convergence est optimale pour une
importante classe d'espaces de Besov anisotropes $B_{2,\infty}^{(\alpha_1,\alpha_2)}$. 
On pr\'esente \'egalement des applications et des simulations.

\end{abstract}

\begin{keyword} adaptive estimation \sep transition density \sep Markov chain \sep model selection
\sep penalized contrast.\\
\emph{2000 MSC} 62M05 \sep 62G05 \sep 62H12 
\end{keyword}


\end{frontmatter}

\section{Introduction}

We consider $(X_i)$ a homogeneous Markov chain. The purpose of this paper is to estimate the 
transition density of such a chain. This quantity allows to comprehend the form of dependence 
between variables and is defined by $\pi(x, y)dy =P(X_{i+1}\in dy | X_i=x)$. It enables also to 
compute other quantities, like $\E[F(X_{i+1})|X_i=x]$  for example. 
As many authors, we choose for this a nonparametric approach. Roussas \cite{Rou69b} first studies an 
estimator of the transition density of a Markov chain. He proves the consistency and the 
asymptotic normality of a kernel estimator for chains satisfying a strong condition known as 
Doeblin's hypothesis. In Bosq \cite{bosq}, an estimator 
by projection is studied in a mixing framework and the consistence is also proved.
 Basu and Sahoo \cite{basu&sahoo} establish a Berry-Essen inequality for a kernel estimator under an assumption 
introduced by Rosenblatt, weaker than the Doeblin's hypothesis.
Athreya and Atuncar \cite{athreya98} improve the result of Roussas since they only need the Harris recurrence
of the Markov chain. 
Other authors are interested in the estimation of the transition density in the non-stationary 
case: Doukhan and Ghind\`es \cite{dou&ghin} bound the integrated risk for any initial distribution. In 
\cite{hernandez}, recursive estimators for a non-stationary Markov chain are described. More 
recently, Clemen\c{c}on \cite{clem} computes the lower bound of the minimax $L^p$ risk
and describes a quotient estimator using wavelets. Lacour \cite{lacour} finds an estimator by projection 
with model selection that reaches the optimal rate of convergence.

All these authors have estimated $\pi$ by observing that $\pi=g/f$ where $g$ is the density of
$(X_i, X_{i+1})$ and $f$ the stationary density. If $\hat{g}$ and $\hat{f}$ are estimators of $g$
and $f$, then an estimator of $\pi$ can be obtained by writing $\hat{\pi}=\hat{g}/\hat{f}$.
But this method has the drawback that the resulting rate of convergence depends on the regularity 
of $f$. And the stationary density $f$ can be less regular than the transition density.

The aim here is to find an estimator $\tilde{\pi}$ of $\pi$ from the observations $X_1,\dots,X_{n+1}$
such that the order of the $L^2$ risk depends only on the regularity of 
$\pi$ and is optimal.

Cl\'emen\c{c}on \cite{clem} introduces an estimation procedure based on an analogy with the regression framework using 
the thresholding of 
wavelets coefficients for regular Markov chains. We propose in this paper an other method based on 
regression, which improves the rate and has the advantage to be really computable. Indeed, this method 
allows to reach the optimal rate of convergence, without the logarithmic loss obtained by Cl\'emen\c{c}on \cite{clem} 
and can be applied to $\beta$-mixing Markov 
chains (the notion of "regular" Markov chains in \cite{clem} is equivalent to $\Phi$-mixing and is then 
a stronger assumption). We use model selection via penalization as described in \cite{BBM} with a new 
contrast inspired by the classical regression contrast. To deal with the dependence
we use auxiliary variables $X_i^*$ as in \cite{viennet97}. But  contrary to most cases 
in such estimation procedure, our penalty does not contain any mixing term and is entirely computable.

In addition, we consider transition densities belonging to anisotropic Besov spaces, i.e. with 
different regularities with respect to the two directions. Our projection spaces (piecewise polynomials, 
trigonometric polynomials or wave\-lets) have different dimensions in the two directions and the procedure 
selects automatically both well fitted dimensions. A lower bound for the rate of convergence on 
anisotropic Besov balls is proved, which shows that our estimation procedure is optimal in 
a minimax sense.

The paper is organized as follows. First, we present the assumptions on the Markov 
chain and on the collections of models. We also give examples of chains and models. 
Section 3 is devoted to estimation procedure and the link with classical regression.
The bound on the empirical risk is established in Section 4 and the $L^2$ control is studied in Section 5.
We compute both upper bound and lower bound for the mean integrated squared error.
In Section 6, some simulation results are given. The proofs are gathered in the last section.


\section{Assumptions}

\subsection{Assumptions on the Markov chain}
We consider an irreducible Markov chain $(X_n)$ taking its values in the real 
line $\mathbb{R}$. We suppose that $(X_n)$ is positive recurrent, i.e. it admits a 
stationary probability measure $\mu$ (for more details, we refer to \cite{m&t}). We assume 
that the distribution $\mu$ has a density $f$ with respect to the Lebesgue measure and that the 
transition kernel $P(x, A)=P(X_{i+1}\in A | X_i=x)$ has also a density, denoted by $\pi$. 
Since the number of observations is finite, $\pi$ is estimated 
on a compact set $A=A_1\times A_2$ only. More precisely, the Markov process is supposed 
to satisfy the following assumptions:

\begin{enumerate}
\item[A1.] $(X_n)$ is irreducible and positive recurrent.
\item[A2.] The distribution of $X_0$ is equal to $\mu$ , thus the chain is (strictly) stationary.
\item[A3.] The transition density $\pi$ is bounded on $A$, i.e.\\
$\|\pi\|_{\infty}:=\sup_{(x,y)\in A}|\pi(x,y)|<\infty$
\item[A4.] The stationary density $f$ verifies $\|f\|_{\infty}:=\sup_{x\in A_1}|f(x)|<\infty$ and  there 
exists a positive real $f_0$ such that, for all $x$ in $A_1$, 
${f(x)\geq f_0}.$
\item[A5.]The chain is geometrically $\beta$-mixing ($\beta_q\leq e^{-\gamma q}$),
or arithmetically $\beta$-mixing ($\beta_q\leq q^{-\gamma}$).
\end{enumerate}
Since $(X_i)$ is a stationary Markov chain, the $\beta$-mixing is very explicit, the mixing coefficients 
can be written: \begin{equation}
\beta_q= \int \|P^q(x,.)-\mu\|_{TV}f(x)dx 
\label{coeffmelange}\end{equation}
where $\|.\|_{TV}$ is the total variation norm (see \cite{doukhan}).

 Notice that we distinguish the sets $A_1$ and $A_2$ in this work because the two directions $x$ and 
$y$ in $\pi(x,y)$ do not play the same role, but in practice $A_1$ and $A_2$ will be equal 
and identical or close to the value domain of the chain. 

\subsection{Examples of chains}

A lot of processes verify the previous assumptions, as (classical or 
more general) autoregressive processes, or diffusions. Here we give a nonexhaustive list of such chains.

\subsubsection{Diffusion processes}
We consider the process $(X_{i\Delta })_{1\leq i \leq n}$ where 
$\Delta>0$ is the observation step and $(X_t)_{t\geq 0}$ is defined by
$$dX_t=b(X_t)dt +\sigma(X_t)dW_t$$
where $W$ is the standard Brownian motion, $b$ is a locally bounded Borel function and
$\sigma$ an uniformly continuous function.
We suppose that the drift function $b$ and the diffusion coefficient $\sigma$ satisfy the 
following conditions, given in \cite{pardouxveretennikov01}(Proposition 1): 
\begin{enumerate}
\item there exists $\lambda_-, \lambda_+$ such that $\quad
    \forall x\neq0, \quad 0<\lambda_-<\sigma^2(x)<\lambda_+$,
\item there exists $M_0\geq 0, \alpha>-1$ and  $r>0$ such that 
      $$\forall |x|\geq M_0, \quad xb(x)\leq -r |x|^{\alpha +1 }.$$
\end{enumerate}

Then, if $X_0$ follows the stationary distribution, 
the discretized process $(X_{i\Delta })_{1\leq i \leq n}$ satisfies Assumptions A1--A5. Note that 
the mixing is geometrical as soon as $\alpha\geq 0$. The continuity of the transition density ensures that 
Assumption A3 holds. Moreover, we can write 
$$f(x)=\frac{1}{M\sigma^2(x)}\exp\left[2\int_0^x\frac{b(u)}{\sigma^2(u)}du\right]$$
with $M$ such that $\int f=1$. Consequently Assumption A4 is verified with 
$\|f\|_\infty\leq \frac{1}{M\lambda_-}\exp\left[\frac{2}{\lambda_-}\sup_{x\in A_1}\int_0^x b(u)du\right]$
and $f_0\geq \frac{1}{M\lambda_+}\exp\left[\frac{2}{\lambda_+}\inf_{x\in A_1}\int_0^xb(u)du\right]$.

\subsubsection{Nonlinear AR(1) processes}
Let us consider the following process
$$X_n=\varphi(X_{n-1})+\varepsilon_{X_{n-1},n}$$
where $\varepsilon_{x,n}$ has a positive density  $l_x$ with respect to the Lebesgue measure, 
which does not depend on $n$. We suppose that the following conditions are verified:
\begin{enumerate}
\item There exist $M>0$ and $\rho<1$ such that, for all $|x|>M$, $|\varphi(x)|<\rho |x|$
and $\sup_{|x|\leq M}|\varphi(x)|<\infty.$
\item There exist $l_0>0$, $l_1>0$ such that $\forall x, y \quad l_0\leq l_x(y)\leq l_1.$
\end{enumerate}
Then Mokkadem \cite{Mokkadem} proves that the chain is Harris recurrent and geometrically ergodic. 
It implies that Assumptions A1 and A5 are satisfied. Moreover 
$\pi(x,y)=l_x(y-\varphi(x))$ and $f(y)=\int f(x)\pi(x,y)dx$ and then Assumptions A3-A4 hold 
with $f_0\geq l_0$ and $\|f\|_\infty\leq \|\pi\|_\infty\leq l_1$.

\subsubsection{ARX(1,1) models}
The nonlinear process ARX(1,1) is defined by
$$X_n=F(X_{n-1},Z_n) + \xi_n$$
where $F$ is bounded and $(\xi_n)$, $(Z_n)$ are independent sequences of i.i.d.
random variables with $\E|\xi_n|<\infty$. We suppose that the distribution of  $Z_n$ has a positive  
density $l$ with respect to the Lebesgue mesure.
Assume that there exist $\rho<1$, a locally bounded and mesurable function 
$h:\mathbb{R}\mapsto\mathbb{R}^+$ such that $\E h(Z_n)<\infty$ and positive constants $M, c$ 
such that 
$$\forall |(u,v)|>M \quad |F(u,v)|<\rho |u|+h(v)-c \text{  and }
\sup_{|x|\leq M}|F(x)|<\infty.$$
Then Doukhan \cite{doukhan} proves (p.102) that  $(X_n)$ is a geometrically $\beta-$mixing process.
We can write 
$$\pi(x,y)=\int l(z)f_\xi(y-F(x,z))dz$$
where $f_\xi $ is the density of $\xi_n$.
So, if we assume furthermore that there exist $a_0,a_1>0$ such that $a_0\leq f_\xi\leq a_1$,
then Assumptions A3-A4 are verified with $f_0\geq a_0$ and $\|f\|_\infty\leq 
\|\pi\|_\infty\leq a_1$.

\subsubsection{ARCH processes}
The model is 
$$X_{n+1}=F(X_n)+ G(X_n)\varepsilon_{n+1}$$
where $F$ and $G$ are continuous functions and for all $x$, $G(x)\neq 0$.
We suppose that the distribution of  $\varepsilon_n$ has a positive density $l$ 
with respect to the Lebesgue measure 
and that there exists $s\geq 1$ such that $\E|\varepsilon_n|^s<\infty$.
The chain $(X_n)$ satisfies Assumptions A1 and A5 if (see \cite{angonze92}):
\begin{equation}
\limsup_{|x|\to\infty} \frac{|F(x)|+|G(x)|(\E|\varepsilon_n|^s)^{1/s}}{|x|}<1.
\label{ARCH}
\end{equation}
In addition, we assume that $\forall x \quad l_0\leq l(x) \leq l_1.$
Then Assumption A3 is verified 
with $\|\pi\|_\infty\leq l_1/\inf_{x \in A_1} G(x)$. And Assumption A4 holds with
$f_0\geq l_0\int f G^{-1}$ and $\|f\|_\infty\leq l_1\int f G^{-1}.$

\subsection{Assumptions on the models}\label{hypmo}
In order to estimate $\pi$, we need to introduce a collection $\{S_m ,m\in \mathcal{M}_n\}$ of 
spaces, that we call models. For each $m=(m_1,m_2)$, $S_m$ is a space of functions with support in $A$
defined from two spaces: $F_{m_1}$ and $H_{m_2}$. $F_{m_1}$ is a subspace of $(L^2\cap L^\infty)(\mathbb{R})$ 
spanned by an orthonormal basis $(\varphi_j^m)_{j\in J_m}$ with $|J_m|=D_{m_1}$ such that, for all $j$, 
the support of $\varphi_j^m$ is included in $A_1$. In the same way $H_{m_2}$ is a subspace of 
$(L^2\cap L^\infty)(\mathbb{R})$ spanned by an orthonormal basis $(\psi_k^m)_{k\in K_m}$ with 
$|K_m|=D_{m_2}$ such that, for all $k$, the support of $\psi_k^m$ is included in $A_2$. 
Here $j$ and $k$ are not necessarily integers, it can be couples of integers as in the case 
of a piecewise polynomial space. Then, we define
$$S_m=F_{m_1}\otimes H_{m_2}=
\{t, \quad t(x,y)=\sum_{j\in J_m}\sum_{k\in K_m} {a}_{j,k}^m \varphi_j^m(x)\psi_k^m(y)\}$$

 The assumptions on the models are the following:

\begin{enumerate}
\item[M1.] For all $m_2$, $D_{m_2}\leq n^{1/3}$
and $\mathcal{D}_n:=\max_{m\in \mathcal{M}_n} D_{m_1}\leq {n}^{1/3}$

\item[M2.] There exist positive reals $\phi_1, \phi_2$ such that, for all $u$ in $F_{m_1}$, 
$\|u\|_\infty^2\leq \phi_1D_{m_1} \int u^2$, and for all $v$ in $H_{m_2}$,
$\sup_{x\in A_2}|v(x)|^2\leq \phi_2D_{m_2}\int v^2. $
By letting $\phi_0=\sqrt{\phi_1\phi_2}$, that leads to
\begin{equation}
\forall t\in S_m \qquad \|t\|_\infty\leq \phi_0\sqrt{D_{m_1}D_{m_2}}\|t\|\label{M2}
\end{equation}
where $\|t\|^2=\int_{\mathbb{R}^2} t^2(x,y)dxdy. $
\item[M3.]$D_{m_1}\leq D_{m_1'}\Rightarrow F_{m_1}\subset F_{m_1'}$ and
$D_{m_2}\leq D_{m_2'}\Rightarrow H_{m_2}\subset H_{m_2'}$
\end{enumerate}

The first assumption guarantees that $\dim S_m=D_{m_1}D_{m_2}\leq n^{2/3}\leq n$ where 
$n$ is the number of observations. The condition M2 implies a useful link
between the $L^2$ norm and the infinite norm.
The third assumption ensures that, for $m$ and $m'$ in $\mathcal{M}_n$,
$S_m+S_{m'}$ is included in a model (since $S_m+S_{m'}\subset S_{m''}$ with 
$D_{m_1''}=\max(D_{m_1},D_{m'_1})$ and $D_{m_2''}=\max(D_{m_2},D_{m'_2})$).
We denote by $\mathcal{S}$ the space with maximal dimension among the $(S_m)_{m\in\mathcal{M}_n}$. 
Thus for all $m$ in $\mathcal{M}_n$, $S_m\subset\mathcal{S}$.

\subsection{Examples of models}\label{ex}
We show here that Assumptions M1--M3 are not too restrictive. Indeed,  
they are verified for the spaces $F_{m_1}$ (and $H_{m_2}$) spanned by the following bases 
(see \cite{BBM}):

\begin{itemize}
\item Trigonometric basis: for $A=[0,1]$,
$<\varphi_0,\dots,\varphi_{m_1-1}>$ with $\varphi_0=\1_{[0,1]}$, 
$\varphi_{2j}(x)=\sqrt{2}$ $\cos(2\pi jx)$ $\1_{[0,1]}(x)$, $\varphi_{2j-1}(x)=\sqrt{2}\sin(2\pi jx)
\1_{[0,1]}(x)$ for $j\geq 1$. For this model $D_{m_1}=m_1$ and $\phi_1=2$ hold.

\item Histogram basis: for $A=[0,1]$,
 $ <\varphi_1,\dots,\varphi_{2^{m_1}}>$ with $\varphi_j
=2^{m_1/2}\1_{[(j-1)/2^{m_1},{j}/2^{m_1}[}  $ for $j=1,\dots,2^{m_1}$. 
Here $D_{m_1}=2^{m_1}$, $\phi_1=1$. 

\item Regular piecewise polynomial basis: for $A=[0,1]$,
polynomials of degree $0,\dots,r$ (where $r$ is fixed) on each interval 
$[(l-1)/2^D,l/2^D[, l=1,\dots,2^D$. In this case, $m_1=(D,r)$, 
$J_m=\{j=(l,d), \quad 1\leq l\leq 2^D, 0\leq d\leq r\}$, $D_{m_1}=(r+1)2^D$. We can put 
$\phi_1=\sqrt{r+1}$.

\item Regular wavelet basis:
$<\Psi_{lk}, l=-1,\dots,m_1, k\in\Lambda(l)>$ where $\Psi_{-1,k}$ points out the translates of 
the father wavelet and $\Psi_{lk}(x)=2^{l/2}\Psi(2^lx-k)$ where $\Psi$ is the mother wavelet.
We assume that the support of the wavelets is included in $A_1$ and that $\Psi_{-1}$ 
belongs to the Sobolev space $W_2^r$.
\end{itemize}


\section{Estimation procedure}\label{sectionest}

\subsection{Definition of the contrast}
To estimate the function $\pi$, we define the contrast \begin{equation}
\gamma_n(t)=\frac{1}{n}\sum_{i=1}^{n}[\int_\mathbb{R} t^2(X_i,y)dy-2t(X_i, X_{i+1})].
\label {contrast}\end{equation}
We choose this contrast because
$$\E\gamma_n(t)=\|t-\pi\|^2_f - \|\pi\|^2_f$$
where $$\|t\|^2_f=\int_{\mathbb{R}^2} t^2(x,y)f(x)dxdy.$$ 
Therefore $\gamma_n(t)$ is the empirical counterpart of the $\|.\|_f$-distance between $t$ and 
$f$ and the minimization of this contrast comes down to minimize $\|t-\pi\|_f$. This contrast is 
new but is actually connected with the one used in regression problems, as we will see in
the next subsection.

We want to estimate $\pi$ by minimizing this contrast on $S_m$.
Let $t(x,y)=\sum_{j\in J_m}\sum_{k\in K_m}$ ${a}_{j,k} \varphi_j^m(x)\psi_k^m(y)$ a function in $S_m$.
Then,
if ${A_m}$ denotes the matrix $({a}_{j,k})_{j\in J_m, k\in K_m}$, 
$$ \forall j_0 \forall k_0\quad 
\frac{\partial \gamma_n(t)}{\partial {a}_{j_0,k_0}}=0\Leftrightarrow G_m{A_m}= Z_m,$$
where $\begin{cases}
G_m=\left(\displaystyle \frac{1}{n}\sum_{i=1}^n\varphi_j^m(X_i)\varphi_l^m(X_i)\right)
_{j, l \in J_m}\\
Z_m=\left(\displaystyle\cfrac{1}{n}\sum_{i=1}^n\varphi_j^m(X_i)\psi_k^m(X_{i+1})\right)
_{j\in J_m, k \in K_m}\\
\end{cases}$\\

Indeed, \begin{equation}
\frac{\partial \gamma_n(t)}{\partial {a}_{j_0,k_0}}=0 \Leftrightarrow 
\sum_{j\in J_m}{a}_{j,k_0}\frac{1}{n}\sum_{i=1}^n\varphi_j^m(X_i)\varphi_{j_0}^m(X_i)
=\cfrac{1}{n}\sum_{i=1}^n\varphi_{j_0}^m(X_i)\psi_{k_0}^m(X_{i+1}) .\label{derivation}
\end{equation}

We can not define a unique minimizer of the contrast $\gamma_n(t)$, since $G_m$ is not 
necessarily invertible. For example, $G_m$ is not invertible if there exists $j_0$ in $J_m$ such that 
there is no observation in the support of $\varphi_{j_0}$ ($G_m$ has a null column). This 
phenomenon happens when localized bases (as histogram bases or piecewise polynomial bases) 
are used. However, the following proposition will enable us to define an estimator:

\begin{prop}\label{existenceminimizer}
$$\forall j_0 \forall k_0\quad 
\frac{\partial \gamma_n(t)}{\partial {a}_{j_0,k_0}}=0 \Leftrightarrow \forall y\quad
(t(X_i,y))_{1\leq i\leq n}=P_W\left(\left(\sum_{k}\psi_k^m(X_{i+1})\psi_k^m(y)\right)_{1\leq i\leq n}\right)$$ 
where $P_W$ denotes the orthogonal projection on $W=\{ (t(X_i,y))_{1\leq i\leq n}, \,t \in S_m\}$.
\end{prop}

Thus the minimization of $\gamma_n(t)$ leads to a unique vector $(\hat{\pi}_m(X_i,y))_{1\leq i\leq n}$ 
defined as 
the projection of $\left(\sum_{k}\psi_k(X_{i+1})\psi_k(y)\right)_{1\leq i\leq n}$ on $W$. The associated 
function $\hat{\pi}_m(.,.)$ is not defined uniquely but we can choose a function $\hat{\pi}_m$ in $S_m$ 
whose values at $(X_i,y)$ are fixed according to Proposition \ref{existenceminimizer}. For the sake 
of simplicity, we denote 
$$\hat\pi_m=\arg\min_{t\in S_m}\gamma_n(t).$$
This underlying function is more a theoretical tool and  the estimator is actually the vector 
$(\hat{\pi}_m(X_i,y))_{1\leq i\leq n}$.
This remark leads to consider the risk defined with the empirical norm 
\begin{equation}
\|t\|_n=\left(\frac{1}{n}\sum_{i=1}^{n}\int_\mathbb{R} t^2(X_i,y)dy\right)^{1/2}.
\label{empiricalnorm}\end{equation}
This norm is the natural distance in this problem and we can notice 
that if $t$ is deterministic with support included in $A_1\times\mathbb{R}$
$$f_0\|t\|^2\leq\E\|t\|_n^2=\|t\|_f^2\leq\|f\|_\infty\|t\|^2$$
and then the mean of this empirical norm is equivalent to the $L^2$ norm $\|.\|$.

\subsection{Link with classical regression}

Let us fix $k$ in $K_m$ and let $$Y_{i,k}=\psi_k^m(X_{i+1}) \qquad\text{ for } i\in\{1,\dots,n\} 
$$
$$t_k(x)=\int t(x,y)\psi_k^m(y)dy \qquad\text{ for all } t \text{ in } L^2(\mathbb{R}^2).$$
Actually, $Y_{i,k}$ and $t_k$ depend on $m$ but we do not mention this for the sake of simplicity.
For the same reason, we denote in this subsection $\psi_k^m$ by $\psi_k$ and $\varphi_j^m$ by 
$\varphi_j$. Then, if $t$ belongs to $S_m$, \begin{eqnarray*}
t(x,y)&=&\sum_{j\in J_m}\sum_{k\in K_m}\left(\int t(x',y')\varphi_j(x')\psi_k(y')dx'dy'\right)
\varphi_j(x)\psi_k(y)\\
&=&\sum_{k\in K_m}\sum_{j\in J_m}\left(\int t_k(x')\varphi_j(x')dx'\right)\varphi_j(x)\psi_k(y)=
\sum_{k\in K_m}t_k(x)\psi_k(y)
\end{eqnarray*}
and then, by replacing this expression of $t$ in $\gamma_n(t)$, we obtain
\begin{eqnarray*}
\hspace{-0.3cm}\gamma_n(t)&=&\frac{1}{n}\sum_{i=1}^{n}[\int \sum_{k,k'}t_k(X_i)t_{k'}(X_i)
\psi_k(y)\psi_{k'}(y)dy
-2\sum_{k}t_k(X_i)\psi_k( X_{i+1})]\\
&=&\frac{1}{n}\sum_{i=1}^{n}\sum_{k\in K_m}[t_k^2(X_i)-2t_k(X_i)Y_{i,k}]
=\frac{1}{n}\sum_{i=1}^{n}\sum_{k\in K_m}[t_k(X_i)-Y_{i,k}]^2-Y_{i,k}^2.
\end{eqnarray*}
Consequently
$$\min_{t\in S_m} \gamma_n(t)=\sum_{k\in K_m}\min_{t_k\in F_{m_1}} 
\frac{1}{n}\sum_{i=1}^{n}[t_k(X_i)-Y_{i,k}]^2-Y_{i,k}^2.$$

We recognize, for all $k$, the least squares contrast, which is used in regression problems.
Here the regression function is $\pi_k=\int \pi(.,y)\psi_k(y)dy$ which verifies
\begin{equation}
Y_{i,k}=\pi_k(X_i) +\varepsilon_{i,k}
\label{mod}\end{equation}
where \begin{equation}\varepsilon_{i,k}=\psi_k(X_{i+1})-\E[\psi_k(X_{i+1})|X_i].\label{defepsilon}
\end{equation}

The estimator $\hat{\pi}_m$ can be written as 
$\sum_{k\in K_m}\hat{\pi}_k(x)\psi_k(y)$ 
where $\hat{\pi}_k$ is the classical least squares estimator for the regression model \eqref{mod} 
(as previously, 
only the vector $(\hat{\pi}_k(X_i))_{1\leq i\leq n}$ is uniquely defined).

This regression model is used in Cl\'emen\c{c}on \cite{clem} to estimate the transition density. 
In the same manner, we could use here the contrast 
$\gamma_n^{(k)}(t)=\frac{1}{n}\sum_{i=1}^{n}[\psi_k(X_{i+1})-t(X_i)]^2$ 
to take advantage of analogy with regression. This method allows to have a good estimation of 
the projection of $\pi$ on some $S_m$ by estimating first each $\pi_k$, but does not provide 
an adaptive method. Model selection requires a more global contrast, as described in~\eqref{contrast}.

\subsection{Definition of the estimator}

We have then an estimator of $\pi$ for all $S_m$.
Let now
$$\hat{m}=\arg\min_{m\in \mathcal{M}_n} \{\gamma_n(\hat{\pi}_m)+\pen(m)\}$$
where $\pen$ is a penalty function to be specified later.
Then we can define $\tilde{\pi}=\hat{\pi}_{\hat{m}}$ and 
compute the empirical mean integrated squared error $\E\|\pi-\tilde{\pi}\|_n^2$
where $\|.\|_n$ is the empirical norm defined in \eqref{empiricalnorm}.


\section{Calculation of the risk}

For a function $h$ and a subspace $S$, let 
$$d(h,S)=\inf_{g\in S}\|h-g\|=\inf_{g\in S}\left(\iint |h(x,y)-g(x,y)|^2dxdy\right)^{1/2}.$$
With an inequality of Talagrand \cite{talagrand1996}, we can prove the following result.
 
\begin{thm} \label{main}
 We consider a Markov chain satisfying Assumptions A1--A5 (with $\gamma>14$ in the case 
of an arithmetical mixing). We consider 
$\tilde{\pi}$ the estimator of the transition density $\pi$ described in Section \ref{sectionest}
with models verifying Assumptions M1--M3 
and the following penalty:\begin{equation}\pen(m)= K_0 \|\pi\|_\infty\frac{D_{m_1}D_{m_2}}{n} 
\label{penalite}\end{equation} where $K_0$ is a numerical constant. Then 
$$\E\|\pi\1_A-\tilde{\pi}\|_n^2\leq C\underset{m\in\mathcal{M}_n}{\inf}\{d^2(\pi\1_A,S_m) 
    +\pen(m)\}+\frac{C'}{n}$$
where $C=\max(5\|f\|_\infty,6)$ and $C'$ is a constant depending on 
$\phi_1, \phi_2, \|\pi\|_\infty, f_0,$ $\|f\|_\infty, \gamma .$ 
\end{thm}

The constant $K_0$ in the penalty is purely numerical (we can choose $K_0=45$).
We observe that the term $\|\pi\|_\infty$ appears in the penalty although it is unknown. 
Nevertheless it can be replaced by any bound of $\|\pi\|_\infty$. Moreover, it is possible 
to use $\|\hat{\pi}\|_\infty$ where $\hat{\pi}$ is some estimator of $\pi$. This method of 
random penalty (specifically with infinite norm) is successfully used in \cite{birge&massart97} 
and \cite{comte2001} for example, and can be applied here even if it means considering 
$\pi$ regular enough. This is proved in appendix.

It is relevant to notice that the penalty term does not contain any mixing term and is then 
entirely computable. It is in fact related to martingale properties of the underlying empirical processes. 
The constant $K_0$ is a fixed universal numerical constant; for practical
purposes, it is adjusted by simulations.

We are now interested in the rate of convergence of the risk. We consider that $\pi$ restricted 
to $A$ belongs to the anisotropic Besov space on $A$ with regularity 
$\boldsymbol{\alpha}=(\alpha_1,\alpha_2)$. Note that if $\pi$ belongs to 
$B_{2,\infty}^{\boldsymbol{\alpha}}(\mathbb{R}^2)$, then $\pi$ restricted 
to $A$ belongs to $B_{2,\infty}^{\boldsymbol{\alpha}}(A)$.
Let us recall the definition of $B_{2,\infty}^{\boldsymbol{\alpha}}(A)$.
Let $e_1$ and $e_2$ be the canonical basis vectors in 
$\mathbb{R}^2$ and for $i=1,2$, $A_{h,i}^r=\{x\in \mathbb{R}^2 ; x, x+he_i, \dots, x+rhe_i \in A\}$.
Next, for $x$ in $A_{h,i}^r$, let
$$\Delta_{h,i}^rg(x)=\sum_{k=0}^r(-1)^{r-k}\binom{r}{k}g(x+khe_i)$$
the $r$th difference operator with step $h$. 
For $t>0$, the directional moduli of smoothness are given by
$$\omega_{r_i,i}(g,t)=\underset{|h|\leq t}{\sup}\left(\int_{A_{h,i}^{r_i}}
|\Delta_{h,i}^{r_i}g(x)|^2dx\right)^{1/2}.$$ We say that $g$ is in the Besov space 
$B_{2,\infty}^{\boldsymbol{\alpha}}(A)$ 
if $$\sup_{t>0}\sum_{i=1}^2t^{-\alpha_i}\omega_{r_i,i}(g,t)<\infty$$ 
for $r_i$ integers larger than $\alpha_i$.
The transition density $\pi$ can thus have different smoothness properties with respect to 
different directions. The procedure described here allows an adaptation of the approximation 
space to each directional regularity. More precisely, if $\alpha_2>\alpha_1$ for example, 
the estimator chooses a space of dimension $D_{m_2}=D_{m_1}^{\alpha_1/\alpha_2}<D_{m_1}$ for 
the second direction, where $\pi$ is more regular. We can thus write the following corollary. 

\begin{cor} \label{coromain}
We suppose that $\pi$ restricted to $A$ belongs to the anisotropic Besov space 
$B_{2,\infty}^{\boldsymbol{\alpha}}(A)$ with regularity $\boldsymbol{\alpha}=(\alpha_1,\alpha_2)$
such that $\alpha_1-2\alpha_2+2\alpha_1\alpha_2>0$ and 
$\alpha_2-2\alpha_1+2\alpha_1\alpha_2>0 $.
We consider the spaces described in Subsection \ref{ex} (with the regularity $r$ of the polynomials
and the wavelets larger than $\alpha_i-1$). Then, 
under the assumptions of Theorem \ref{main}, 
 $$\E\|\pi\1_A-\tilde{\pi}\|_n^2=O(n^{-\frac{2\bar\alpha}{2\bar\alpha+2}}).$$
where $\bar\alpha$ is the harmonic mean of $\alpha_1$ and $\alpha_2$.
\end{cor}

The harmonic mean of $\alpha_1$ and $\alpha_2$ is the real $\bar\alpha$ such that
$2/\bar\alpha=1/\alpha_1+1/\alpha_2.$
Note that the condition $\alpha_1-2\alpha_2+2\alpha_1\alpha_2>0 $ 
is ensured as soon as $\alpha_1\geq 1$ and the condition $\alpha_2-2\alpha_1+2\alpha_1\alpha_2>0 $ 
as soon as $\alpha_2\geq 1$. 

Thus we obtain the rate of convergence $n^{-\frac{2\bar\alpha}{2\bar\alpha+2}}$, which is optimal 
in the minimax sense (see Section 5.3 for the lower bound).


\section{$L^2$ control}

\subsection{Estimation procedure}\label{estL2}

Although  the empirical norm is the more natural in this problem, we are interested in a $L^2$
control of the risk. For this, the estimation procedure must be modified. 
We truncate the previous estimator in the following way :
\begin{equation}
\tilde{\pi}^*=\begin{cases}
\tilde{\pi}& \text{ if }\|\tilde{\pi}\|\leq k_n\\
0& \text{ else }\end{cases}
\label{est}\end{equation}
with $k_n= n^{2/3}$.

\subsection{Calculation of the $L^2$ risk}

We obtain in this framework a result similar to Theorem \ref{main}.

\begin{thm} \label{main2}
We consider a Markov chain satisfying Assumptions A1--A5 (with $\gamma>20$ in the case of 
an arithmetical mixing). We consider 
$\tilde{\pi}^*$ the estimator of the transition density $\pi$ described in Section \ref{estL2}.
Then 
$$\E\|\tilde{\pi}^*-\pi\1_A\|^2\leq C\underset{m\in\mathcal{M}_n}{\inf}\{d^2(\pi\1_A,S_m) 
    +\pen(m)\}+\frac{C'}{n}.$$
where $C=\max(36f_0^{-1}\|f\|_\infty+2,36f_0^{-1})$ and 
$C'$ is a constant depending on $\phi_1, \phi_2, \|\pi\|_\infty, \|\pi\|,f_0, \|f\|_\infty, \gamma.$
\end{thm}

If $\pi$ is regular, we can state the following corollary:

\begin{cor}
We suppose that the restriction of $\pi$ to $A$ belongs to the aniso\-tropic Besov space 
$B_{2,\infty}^{\boldsymbol{\alpha}}(A)$ with regularity $\boldsymbol{\alpha}=(\alpha_1,\alpha_2)$
such that $\alpha_1-2\alpha_2+2\alpha_1\alpha_2>0 $ and $\alpha_2-2\alpha_1+2\alpha_1\alpha_2>0$.
We consider the spaces described in Subsection \ref{ex} (with the regularity $r$ of the polynomials
and the wavelets larger than $\alpha_i-1$). Then, under the assumptions of Theorem \ref{main2}, 
$$\E\|\pi\1_A-\tilde{\pi}^*\|^2=O(n^{-\frac{2\bar\alpha}{2\bar\alpha+2}}).$$
where $\bar\alpha$ is the harmonic mean of $\alpha_1$ and $\alpha_2$.
\end{cor}

The same rate of convergence is then achieved with the $L^2$ norm instead of the empirical norm. And the procedure allows to 
adapt automatically the two dimensions of the projection spaces to the regularities $\alpha_1$ 
and $\alpha_2$ of the transition density $\pi$.
If $\alpha_1=1$ we recognize the rate $n^{-\frac{\alpha_2}{3\alpha_2+1}}$ established by Birg\'e \cite{birge83} 
with metrical arguments. 
The optimality is proved in the following subsection.

If $\alpha_1=\alpha_2=\alpha$ ("classical" Besov space), then $\bar\alpha=\alpha$ and 
our result is thus an improvement of the one of Cl\'emen\c{c}on \cite{clem}, whose 
procedure achieves only the rate $(\log(n)/n)^{\frac{2\alpha}{2\alpha+2}}$ and allows 
to use only wavelets.
We can observe that in this case, the condition $\alpha_1-2\alpha_2+2\alpha_1\alpha_2>0 $ 
is equivalent to $\alpha>1/2$ and so is verified if  the function $\pi$ is regular enough.\\
Actually, in the case $\alpha_1=\alpha_2$, an estimation with isotropic spaces ($D_{m_1}=D_{m_2}$)
is preferable. Indeed, in this framework, the models are nested and so we can consider spaces 
with larger dimension ($D_m^2\leq n$ instead of $D_m^2\leq n^{2/3}$). Then Corollary \ref{coromain}
is valid whatever $\alpha>0$. Moreover, for the arithmetic mixing, assumption $\gamma>6$ is sufficient.

\subsection{Lower bound}\label{lower}

We denote by $\|.\|_A$ the norm in $L^2(A)$, i.e. $\|g\|_{A}=\left(\int_A|g|^2\right)^{1/2}$.
We set
\begin{eqnarray*}
  \mathcal{B}=\{\pi \text{ transition density on }\mathbb{R} \text{ of a positive recurrent } \\
\text{ Markov chain such that } \|\pi\|_{B_{2,\infty}^{\boldsymbol{\alpha}}(A)}\leq L\}
\end{eqnarray*}
and $\E_\pi$ the expectation corresponding to the distribution of $X_1, \dots, X_n$ 
if the true transition density of the Markov chain is $\pi$ and the initial distribution is
the stationary distribution.

\begin{thm}
There exists a positive constant $C$ such that, if $n$ is large enough,
$$\inf_{\hat{\pi}_n}\sup_{\pi \in \mathcal{B}} \E_\pi\|\hat{\pi}_n-\pi\|_{A}^2\geq Cn^{-\frac{2\bar\alpha}{2\bar\alpha+2}}$$
where the infimum is taken over all estimators $\hat{\pi}_n$ of $\pi$ based on the observations 
$ X_1,\dots, X_{n+1}$. 
\label{lowerbound}\end{thm}

So the lower bound in \cite{clem} is generalized for the case $\alpha_1\neq\alpha_2$. 
It shows that our procedure reaches the optimal minimax rate, 
whatever the regularity of $\pi$, without needing to know $\boldsymbol{\alpha}$. 

\section{Simulations}

To evaluate the performance of our method, we simulate a Markov chain with a known transition density
and then we estimate this density and compare the two functions for different values of $n$.
The estimation procedure is easy, we can decompose it in some steps:

\begin{itemize}
\item find the coefficients matrix $A_m$ for each $m=(m_1,m_2)$
\item compute $\gamma_n(\hat\pi_m)=\mathrm{Tr}(\transpose{A_m}G_mA_m-2\transpose{Z_m}A_m)$
\item find $\hat{m}$ such that $\gamma_n(\hat\pi_m)+\pen(m)$ is minimum
\item compute $\hat\pi_{\hat{m}}$
\end{itemize}

For the first step, we use two different kinds of bases : the histogram bases and the trigonometric bases, 
as described in subsection \ref{ex}. We renormalize these bases so that they are defined on the estimation 
domain $A$ instead of $[0,1]^2$.
For the third step, we choose $\pen(m)=0.5\cfrac{D_{m_1}D_{m_2}}{n}$.

We consider three Markov chains: \\
$\bullet$ An autoregressive process defined by $X_{n+1}=aX_n+b+\varepsilon_{n+1}$, where 
the $\varepsilon_{n}$ are i.i.d. centered Gaussian random variables with variance $\sigma^2$.
The stationary distribution of this process is a Gaussian with mean $b/(1-a)$ and with variance
$\sigma^2/(1-a^2)$. The transition density is $\pi(x,y)=\varphi(y-ax-b)$ where 
$\varphi(z)=1/(\sigma\sqrt{2\pi}).\exp(-z^2/2\sigma^2)$ is the density of a standard Gaussian.
Here we choose $a=0.5$, $b=3$, $\sigma=1$ and we note this process AR(1). It is estimated
on $[4,8]^2$.\\
$\bullet$ A discrete radial Ornstein-Uhlenbeck process, i.e. the Euclidean norm of a vector 
$(\xi^1,\xi^2,\xi^3)$ whose components are i.i.d. processes satisfying, for $j=1,2,3$,
$\xi^j_{n+1}=a\xi^j_n+\beta\varepsilon_n^j$ where $\varepsilon_{n}^j$ are i.i.d. standard Gaussian.
This process is studied in detail in \cite{chaleyatgenon06}. Its transition density is 
$$\pi(x,y)=\1_{y>0}\exp(-\frac{y^2+a^2x^2}{2\beta^2})I_{1/2}(\frac{axy}{\beta^2})\frac{y}{\beta^2}
\sqrt{\frac{y}{ax}}$$ where $I_{1/2}$ is the Bessel function with index $1/2$.
The stationary density of this chain is $f(x)=\1_{x>0}\exp\{-x^2/2\rho^2\}{2x^2}/(\rho^3\sqrt{2\pi})$
with $\rho^2=\beta^2/(1-a^2)$.
We choose $a=0.5$, $\beta=3$ and we denote this process by $\sqrt{\text{CIR}}$ since it is the 
square root of a Cox-Ingersoll-Ross process. The estimation domain is $[2,10]^2$.\\
$\bullet$ An ARCH process defined by $X_{n+1}=\sin(X_n)+(\cos(X_n)+3)\varepsilon_{n+1}$ where
the $\varepsilon_{n+1}$ are i.i.d. standard Gaussian. We verify that the condition \eqref{ARCH} is satisfied.
Here the transition density is $$\pi(x,y)=\varphi\left(\frac{y-\sin(x)}{\cos(x)+3}\right)
\frac{1}{\cos(x)+3}$$ and we estimate this chain on $[-6,6]^2$.

\begin{figure}[!h]\begin{center}
\includegraphics[scale=0.65]{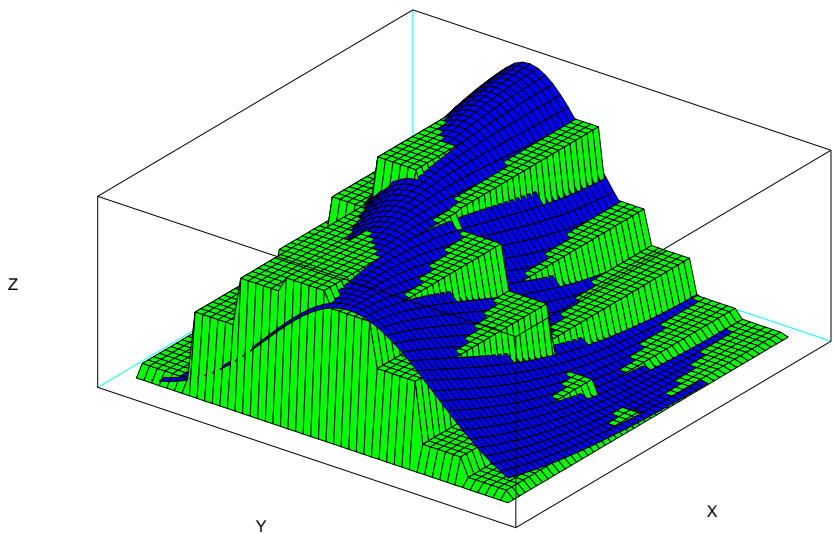}
\caption{Estimator (light surface) and true fonction (dark surface) for a $\sqrt{\text{CIR}}$ process 
estimated with a histogram basis, $n=1000$.}
\end{center}\label{cx}\end{figure}

\begin{figure}[!h]
\begin{tabular}{cc}
\hspace{-0.7cm}\includegraphics[scale=0.7]{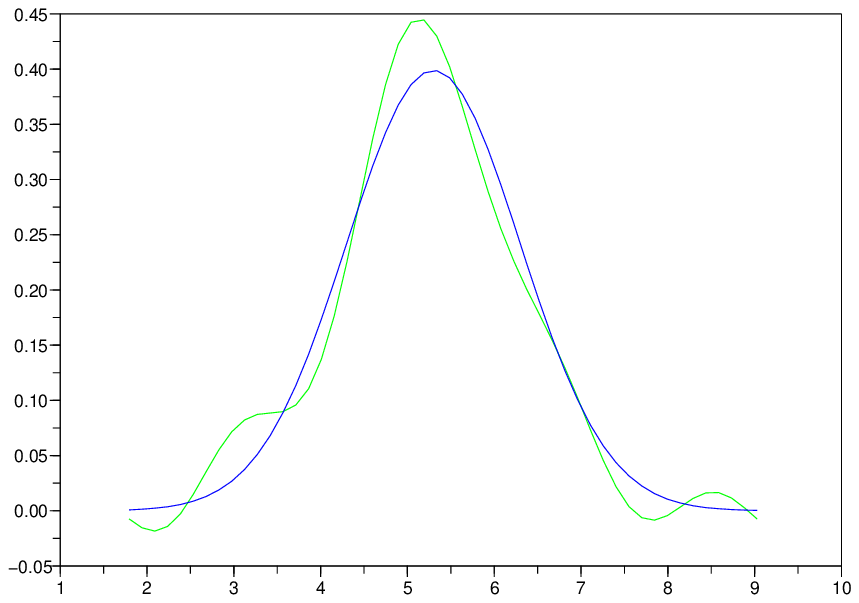} &
\includegraphics[scale=0.7]{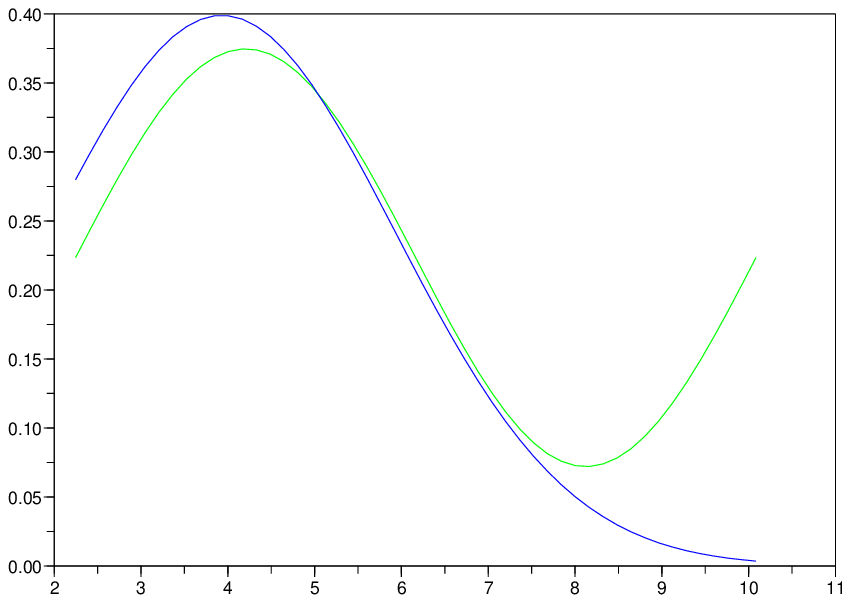}\\
${ x=4.6}$  & ${ y=5}$
\end{tabular}
\caption{Sections for AR(1) process estimated with a trigonometric basis, $n=1000$, dark line: true function,
light line: estimator.}
\label{ar}\end{figure}

We can illustrate the results by some figures. Figure 1 shows the surface $z=\pi(x,y)$
and the estimated surface $z=\tilde{\pi}(x,y)$. We use a histogram basis and we see that the 
procedure chooses different dimensions on the abscissa and on the ordinate since the estimator is constant 
on rectangles instead of squares. Figure \ref{ar} presents sections of this kind of surfaces for the AR(1)
process estimated with trigonometric bases. We can see the curves $z=\pi(4.6,y)$ versus 
$z=\tilde{\pi}(4.6,y)$ and the curves $z=\pi(x,5)$ versus $z=\tilde{\pi}(x,5)$. The second section shows 
that it may exist some edge effects due to the mixed control of the two directions.

For more precise results, empirical risk and $L^2$ risk are given respectively 
in Table \ref{erreursemp} and Table \ref{erreursL2}.

\begin{table}[!h]
\begin{center}
\begin{tabular}{|l|l|l|l|l|l|c|}
\hline
\backslashbox{law}{$n$} 
& 50 & 100 & 250 & 500 & 1000 & basis \\
\hline \hline
AR(1) & 0.067 & 0.055 & 0.043 & 0.038 & 0.033 & H\\
& 0.096 & 0.081 & 0.063 & 0.054 & 0.045 & T \\
\hline
$\sqrt{\text{CIR}}$ & 0.026 & 0.023 & 0.019 & 0.016 & 0.014 & H\\
 & 0.019 & 0.015 & 0.009 & 0.007 & 0.006 & T \\
\hline
ARCH & 0.031 & 0.027 & 0.016 & 0.015 & 0.014 &  H\\
 & 0.020 & 0.012 & 0.008 & 0.007 & 0.007 & T \\
\hline
\end{tabular}
\caption{Empirical risk $\E\|\pi-\tilde{\pi}\|_n^2$ for simulated data with 
$\pen(m)=0.5D_{m_1}D_{m_2}/n$, averaged over $N=200$ samples. H: histogram basis, 
T: trigonometric basis.}\label{erreursemp}
\end{center}
\end{table}

\begin{table}[!h]
\begin{center}
\begin{tabular}{|l|l|l|l|l|l|c|}
\hline
\backslashbox{law}{$n$} 
& 50 & 100 & 250 & 500 & 1000 &  basis \\
\hline \hline
AR(1) & 0.242 & 0.189 & 0.132 & 0.109 & 0.085  & H \\
 & 0.438 & 0.357 & 0.253 & 0.213 & 0.180  & T\\
\hline
$\sqrt{\text{CIR}}$  & 0.152 & 0.130 & 0.094 & 0.066 & 0.054  & H\\
 & 0.152 & 0.123 & 0.072 & 0.052 & 0.046   & T\\
\hline
ARCH & 0.367 & 0.303 & 0.168 & 0.156 & 0.144   & H \\
 & 0.249 & 0.137 & 0.096 & 0.092 & 0.090  & T\\
\hline
\end{tabular}
\caption{$L^2$ risk $\E\|\pi-\tilde{\pi}^*\|^2$ for simulated data with $\pen(m)=0.5D_{m_1}D_{m_2}/n$, 
averaged over $N=200$ samples. H: histogram basis, T: trigonometric basis.}\label{erreursL2}
\end{center}
\end{table}

\begin{table}[!h]
\begin{center}
\begin{tabular}{|l|l|l|l|l|l|c|}
\hline
\backslashbox{law}{$n$} 
& 50 & 100 & 250 & 500 & 1000 &  basis \\
\hline \hline
AR(1) & 0.052 & 0.038 & 0.026 & 0.020 & 0.015 & H \\
 & 0.081 & 0.069 & 0.046 & 0.037 & 0.031 & T\\
\hline
$\sqrt{\text{CIR}}$  & 0.016 & 0.014 & 0.010 & 0.006 & 0.004 & H\\
 & 0.018 & 0.012 & 0.008 & 0.006 & 0.004 & T\\
\hline
\end{tabular}
\caption{$L^2(f(x)dxdy)$ risk $\E\|\pi-\tilde{\pi}^*\|_f^2$ for simulated data with $\pen(m)=0.5D_{m_1}D_{m_2}/n$, 
averaged over $N=200$ samples. H: histogram basis, T: trigonometric basis.}\label{erreursLf}
\end{center}
\end{table}

We observe that the results are better when we consider the empirical norm. It was expectable, given
that this norm is adapted to the studied problem. Actually the better norm to evaluate 
the distance between $\pi$ and its estimator is the norm $\|.\|_f$. Table \ref{erreursLf} shows that 
the errors in this case are very satisfactory. 

So the results are roughly good but we can not pretend that a basis among the others gives better 
results. We can then imagine a mixed strategy, i.e. a procedure which uses several kinds of bases 
and which can choose the best basis. 
These techniques are successfully used in a regression framework by Comte and Rozenholc 
 \cite{comterozenholc2002}, \cite{comterozenholc2004}.


\section{Proofs}

\subsection{Proof of Proposition \ref{existenceminimizer}}
Equality \eqref{derivation} yields, by multiplying by $\psi_{k_0}^m(y)$,
\begin{eqnarray*}
\sum_{j\in J_m}{a}_{j,k_0}\sum_{i=1}^n\varphi_j^m(X_i)\psi_{k_0}^m(y)\varphi_{j_0}^m(X_i)
=\sum_{i=1}^n\varphi_{j_0}^m(X_i)\psi_{k_0}^m(X_{i+1})\psi_{k_0}^m(y).
\end{eqnarray*}
Then, we sum over $k_0$ in $K_m$:
\begin{eqnarray*}
\sum_{i=1}^n t(X_i,y)\varphi_{j_0}^m(X_i)
=\sum_{i=1}^n\sum_{k_0\in K_m}\psi_{k_0}^m(X_{i+1})\psi_{k_0}^m(y)\varphi_{j_0}^m(X_i).
\end{eqnarray*}
If we multiply this equality by $a'_{j_0,k}\psi_{k}^m(y)$ and 
if we sum over $k\in K_m$ and $j_0 \in J_m$, we obtain
\begin{eqnarray*}
\sum_{i=1}^n [t(X_i,y)-\sum_{k_0\in K_m}\psi_{k_0}^m(X_{i+1})\psi_{k_0}^m(y)]
\sum_{k\in K_m}\sum_{j_0 \in J_m} a'_{j_0,k}\varphi_{j_0}^m(X_i)\psi_{k}^m(y)=0\\
\text{ i.e.  }\hspace{3cm}\sum_{i=1}^n [t(X_i,y)-\sum_{k_0\in K_m}\psi_{k_0}^m(X_{i+1})\psi_{k_0}^m(y)]
u(X_i,y)=0
\end{eqnarray*}
for all $u$ in $S_m$.
So the vector $(t(X_i,y)-\sum_{k\in K_m}\psi_{k}^m(X_{i+1})\psi_k^m(y))
      _{1\leq i\leq n}$ is orthogonal to each vector in $W$.
Since $t(X_i,y)$ belongs to $W$, the proposition is proved.

\subsection{Proof of Theorem \ref{main}}

For $\rho$ a real larger than 1, let $$\displaystyle \Omega_\rho
=\{\forall t\in \mathcal{S} \quad \|t\|_f^2\leq\rho\|t\|_n^2\}$$

In the case of an arithmetical mixing, since $\gamma>14$, there exists a real $c$ such that 
$$\begin{cases}
0<c<\cfrac{1}{6}\\
\gamma c>\cfrac{7}{3}
\end{cases}$$
We set in this case $q_n=\frac{1}{2}\lfloor n^c \rfloor.$  
In the case of a geometrical mixing, we set $q_n=\frac{1}{2}\lfloor c\ln(n) \rfloor$ 
where $c$ is a real larger than $7/3\gamma$.

For the sake of simplicity, we suppose that $n=4p_nq_n$, with $p_n$ an integer.
Let for $i=1,\dots,n/2$, $U_i=(X_{2i-1},X_{2i}).$\\
Let $\begin{cases}
A_l=& (U_{2lq_n+1},..., U_{(2l+1)q_n})\qquad l=0,\dots,p_n-1,\\
B_l=& (U_{(2l+1)q_n+1},..., U_{(2l+2)q_n})\qquad l=0,\dots,p_n-1.\\
\end{cases}$\\
We use now the mixing assumption A5. As in Viennet \cite{viennet97}
we can build a sequence $(A_l^*)$ such that 
$$\begin{cases}
A_l \text{ and } A_l^* \text{ have the same distribution},\\
A_l^* \text { and } A_{l'}^* \text{ are independent if }l\neq l',\\
P(A_l\neq A_{l}^*)\leq \beta_{2q_n}.
\end{cases}$$
In the same way, we build $(B_l^*)$ and we define for any $l\in\{0,\dots, p_n-1\}$, \\
$A_l^*= (U^*_{2lq_n+1},..., U^*_{(2l+1)q_n})$, $B_l^*= (U^*_{(2l+1)q_n+1},..., U^*_{(2l+2)q_n})$
so that the sequence $(U_1^*,\dots, U_{n/2}^*)$ and then the sequence $(X_1^*,\dots, X_{n}^*)$
are well defined.\\
Let now $V_i=(X_{2i},X_{2i+1})$ for $i=1,\dots,n/2$ and \\$\begin{cases}
C_l=& (V_{2lq_n+1},..., V_{(2l+1)q_n})\qquad l=0,\dots,p_n-1,\\
D_l=& (V_{(2l+1)q_n+1},..., V_{(2l+2)q_n})\qquad l=0,\dots,p_n-1.\\
\end{cases}$\\
We can build $(V_1^{**},\dots, V_{n/2}^{**})$ and then $(X_2^{**},\dots, X_{n+1}^{**})$ such that 
$$\begin{cases}
C_l \text{ and } C_l^{**} \text{ have the same distribution},\\
C_l^{**} \text { and } C_{l'}^{**} \text{ are independent if }l\neq l',\\
P(C_l\neq C_{l}^{**})\leq \beta_{2q_n}.
\end{cases}$$

We put $X_{n+1}^*=X_{n+1}$ and $X_1^{**}=X_1$. Now let $$\Omega^*=\{\forall i \quad X_i=X_i^*=X_i^{**}\}
\quad \text{  and  } \quad \Omega_\rho^*=\Omega_\rho\cap\Omega^*$$
We denote by $\pi_m$ the orthogonal projection of $\pi$ on $S_m$.
Now, 
\begin{eqnarray}\label{firstsecond}
\E\|\tilde{\pi}-\pi\1_{A}\|_n^2&=&\E\left(\|\tilde{\pi}-\pi\1_{A}\|_n^2\1_{\Omega_{\rho}^*}\right)
+\E\left(\|\tilde{\pi}-\pi\1_{A}\|_n^2\1_{\Omega_{\rho}^{*c}}\right)
\end{eqnarray}

To bound the first term, we observe that for all $s,t$
$$\gamma_n(t)-\gamma_n(s)=\|t-\pi\|_n^2-\|s-\pi\|_n^2 -2 Z_n(t-s)$$
where $\displaystyle Z_n(t)=\frac{1}{n}\sum_{i=1}^{n}\left\{t(X_i,X_{i+1})-
\int_\mathbb{R} t(X_i,y)\pi(X_i,y)dy\right\}.$\\
Since $\|t-\pi\|_n^2=\|t-\pi\1_{A}\|_n^2+\|\pi\1_{A^c}\|_n^2$, we can write
$$\gamma_n(t)-\gamma_n(s)=\|t-\pi\1_{A}\|_n^2-\|s-\pi\1_{A}\|_n^2 -2 Z_n(t-s).$$
The definition of $\hat{m}$ gives, for some fixed $m \in \mathcal{M}_n$,
$$\gamma_n(\tilde{\pi})+\pen(\hat{m})\leq \gamma_n(\pi_m)+\pen(m)$$
And then \begin{eqnarray*}
\|\tilde{\pi}-\pi\1_{A}\|_n^2\leq \|\pi_m-\pi\1_{A}\|_n^2
  +2Z_n(\tilde{\pi}-\pi_m)+\pen(m)-\pen(\hat{m})\\
\leq\|\pi_m-\pi\1_{A}\|_n^2
  +2\|\tilde{\pi}-\pi_m\|_f\sup_{t\in B_f(\hat{m})}Z_n(t)+\pen(m)-\pen(\hat{m})
\end{eqnarray*}
where, for all $m'$, $B_f(m')=\{t\in S_m+S_{m'},\quad \|t\|_f=1\}$.
Let $ \theta$ a real larger than $2\rho$   and $p(.,.)$ a function such that 
$\theta p(m,m')\leq \pen(m)+\pen(m')$. Then
\begin{eqnarray}
\|\tilde{\pi}-\pi\1_{A}\|_n^2\1_{\Omega_{\rho}^*}&\leq& \|\pi_m-\pi\1_{A}\|_n^2
  +\frac{1}{\theta}\|\tilde{\pi}-\pi_m\|_f^2\1_{\Omega_{\rho}^*}+2\pen(m)\nonumber\\ &&
 + \theta\sum_{m'\in \mathcal{M}_n}\left[\sup_{t\in B_f(m')}Z_n^2(t)-p(m,m')\right]_+\1_{\Omega_{\rho }^*} 
\label{in1}\end{eqnarray}
But 
$\|\tilde{\pi}-\pi_m\|_f^2\1_{\Omega_{\rho}^*}\leq\rho\|\tilde{\pi}-\pi_m\|_n^2
\1_{\Omega_{\rho}^*}
\leq2\rho\|\tilde{\pi}-\pi\1_{A}\|_n^2\1_{\Omega_{\rho}^*}+2\rho\|\pi\1_{A}-\pi_m\|_n^2.$

Then, inequality \eqref{in1} becomes 
\begin{eqnarray}\label{ineg}
\|\tilde{\pi}-\pi\1_{A}\|_n^2\1_{\Omega_{\rho}^*}\left(1-\frac{2\rho}{\theta}\right)&\leq&
\left(1+\frac{2\rho}{\theta}\right)\|\pi_m-\pi\1_{A}\|_n^2+2 \pen(m) \nonumber\\ &&
+\theta\sum_{m'\in \mathcal{M}_n}\left[\sup_{t\in B_f(m')}Z_n^2(t)-p(m,m')\right]_+\1_{\Omega_{\rho}^*}\nonumber\\
\text{so }\hspace{0.5cm} \E\left(\|\tilde{\pi}-\pi\1_{A}\|_n^2\1_{\Omega_{\rho }^*}\right)&\leq&
\frac{\theta+2\rho}{\theta-2\rho}
\E\|\pi\1_{A}-\pi_m\|_n^2+\frac{2\theta}{\theta-2\rho}\pen(m)\nonumber  \\ &&
\hspace{-3cm}+\frac{\theta^2}{\theta-2\rho }\sum_{m'\in \mathcal{M}_n}
\E\left(\left[\sup_{t\in B_f(m')}Z_n^2(t)-p(m,m')\right]_+\1_{\Omega_{\rho }^*}\right)\end{eqnarray}

We now use the following proposition:

\begin{prop}\label{talagrand}
Let $p(m,m')=10\|\pi\|_\infty\cfrac{D(m,m')}{n}$ where $D(m,m')$ denotes the dimension of 
$S_m+S_{m'}$. Then, under the assumptions of Theorem \ref{main},
there exists a constant $C_1$ such that  \begin{equation}\sum_{m'\in \mathcal{M}_n}
\E\left(\left[\sup_{t\in B_f(m')}Z_n^2(t)-p(m,m')\right]_+\1_{\Omega^*}\right)\leq \frac{C_1}{n}.
\label{inegtalagrand}\end{equation}
\end{prop}

Then, with $\theta= 3\rho$, inequalities \eqref{ineg} and \eqref{inegtalagrand} yield \begin{equation}
\E\left(\|\tilde{\pi}-\pi\1_{A}\|_n^2\1_{\Omega_{\rho }^*}\right)\leq 5\|f\|_\infty
\|\pi_m-\pi\1_{A}\|^2+6\pen(m)+\frac{9\rho C_1}{n}
\label{firstterm}\end{equation}

The penalty term $\pen(m)$ has to verify $\pen(m)+\pen(m')\geq30\rho\|\pi\|_\infty$
$\cfrac{D(m,m')}{n}$ i.e.  $ 30\rho\|\pi\|_\infty {\rm dim}(S_m+S_{m'})\leq \pen(m)+\pen(m')$
We choose $\rho=3/2$ and so $\pen(m)=45\|\pi\|_\infty\cfrac{D_{m_1}D_{m_2}}{n}$.

To bound the second term in \eqref{firstsecond}, we recall (see Section \ref{sectionest}) that 
$(\hat{\pi}_{\hat{m}}(X_i,y))_{1\leq i\leq n}$ is  
the orthogonal projection of $(\sum_{k}\psi_k(X_{i+1})\psi_k(y))_{1\leq i\leq n}$ on
$$W=\{ (t(X_i,y))_{1\leq i\leq n}, \quad t \in S_{\hat{m}}\}$$
where $\psi_k=\psi_k^{\hat{m}}$.
Thus, since $P_W$ denotes the orthogonal projection on $W$, using \eqref{mod}-\eqref{defepsilon}
\begin{eqnarray*}
(\hat{\pi}_{\hat{m}}(X_i,y))_{1\leq i\leq n}
&=&{P}_W((\sum_{k}\psi_k(X_{i+1})\psi_k(y))_{1\leq i\leq n})\\
&=&{P}_W((\sum_{k}\pi_k(X_{i})\psi_k(y))_{1\leq i\leq n})
+{P}_W((\sum_{k}\varepsilon_{i,k}\psi_k(y))_{1\leq i\leq n})\\
&=&{P}_W(\pi\1_A(X_{i},y))_{1\leq i\leq n})
+{P}_W((\sum_{k}\varepsilon_{i,k}\psi_k(y))_{1\leq i\leq n})
\end{eqnarray*}

We denote by $\|.\|_{\mathbb{R}^n}$ the Euclidean norm in $\mathbb{R}^n$, by $X$ the vector 
$(X_{i})_{1\leq i\leq n}$ and by $\varepsilon_k$ the vector $(\varepsilon_{i,k})_{1\leq i\leq n}$.
Thus 
\begin{eqnarray*}
&&\hspace{-0.5cm}\|\pi\1_A-\hat{\pi}_{\hat{m}}\|_n^2=\frac{1}{n}\int \|\pi\1_A(X,y)-{P}_W(\pi\1_A(X,y))
-{P}_W(\sum_{k}\varepsilon_k\psi_k(y))\|_{\mathbb{R}^n}^2dy\\
&&=\frac{1}{n}\int \|\pi\1_A(X,y)-{P}_W(\pi\1_A(X,y))\|_{\mathbb{R}^n}^2dy
+\frac{1}{n}\int \|{P}_W(\sum_{k}\varepsilon_k\psi_k(y))\|_{\mathbb{R}^n}^2dy\\
&&\leq\frac{1}{n}\int \|\pi\1_A(X,y)\|_{\mathbb{R}^n}^2dy
+\frac{1}{n}\int\|\sum_{k}\varepsilon_k\psi_k(y)\|_{\mathbb{R}^n}^2dy\\
&&\leq \frac{1}{n}\sum_{i=1}^n\|\pi\|_\infty\int\pi(X_i,y)dy+
\frac{1}{n}\sum_{i=1}^n\int[\sum_{k}\varepsilon_{i,k}\psi_k(y)]^2dy\\
&&\leq \|\pi\|_\infty +\frac{1}{n}\sum_{i=1}^n\sum_{k}\varepsilon_{i,k}^2.
\end{eqnarray*}
But Assumption M2 implies 
$\|\sum_{k \in K_{\hat{m}}} \psi_{k}^2\|_\infty\leq \phi_2 D_{\hat{m}_2}$.
So, using \eqref{defepsilon},\begin{eqnarray*}
\varepsilon_{i,k}^2 &\leq &2\psi_k^2(X_{i+1})+2\E[\psi_k(X_{i+1})|X_i]^2 \\
 {\rm and }\quad\sum_{k}\varepsilon_{i,k}^2 &\leq& 2\sum_{k}\psi_k^2(X_{i+1})+2\E[\sum_{k}\psi_k^2(X_{i+1})|X_i]
\leq 4 \phi_2 D_{\hat{m}_2}
\end{eqnarray*}

Thus we obtain \begin{equation}\label{majdepi-pichapeau}
\|\pi\1_A-\hat{\pi}_{\hat{m}}\|_n^2\leq\|\pi\|_\infty+4\phi_2D_{\hat{m}_2}\leq\|\pi\|_\infty+4\phi_2n^{1/3}
\end{equation}
and, by taking the expectation, $\E\left(\|\pi\1_A-\hat{\pi}_{\hat{m}}\|_n^2\1_{\Omega_{\rho}^{*c}}\right)\leq 
(\|\pi\|_\infty+4\phi_2 n^{1/3})P({\Omega_{\rho}^{*c}})$.

We now remark that 
$P(\Omega_\rho^{*c})=P(\Omega^{*c})+P(\Omega_\rho^c\cap\Omega^*).$
In the geometric case $\beta_{2q_n}\leq e^{-\gamma c\ln(n)}\leq n^{-\gamma c}$ and in the other case
$\beta_{2q_n}\leq (2q_n)^{-\gamma}\leq n^{-\gamma c}$. 
Then $$P(\Omega^{*c})\leq 4p_n\beta_{2q_n}\leq n^{1-c\gamma}.$$ 
But we have choosed $c$ such that $c\gamma>{7}/{3}$ 
and so $P(\Omega^{*c})\leq n^{-4/3}.$
Now we will use the following proposition:

\begin{prop}\label{omegac} Let $\rho>1$. Then, under the assumptions of Theorem \ref{main}
or Theorem \ref{main2}, 
there exists $C_2>0$ such that $P(\Omega_{\rho}^{c}\cap\Omega^*)\leq \cfrac{C_2}{n^{7/3}}.$
\end{prop}

This proposition implies that
$\E\left(\|\pi\1_A-\hat{\pi}_{\hat{m}}\|_n^2\1_{\Omega_{\rho}^{*c}}\right)\leq \cfrac{C_3}{n}$.

Now we use \eqref{firstterm} and we observe that 
this inequality holds for all $m$ in $\mathcal{M}_n$, so
$$\E\|\tilde{\pi}-\pi\1_A\|_n^2\leq C\inf_{m\in \mathcal{M}_n}(\|\pi\1_A-\pi_m\|^2+\pen(m))+\frac{C_4}{n}$$
with $C=\max(5\|f\|_\infty,6)$.

\subsection{Proof of Corollary \ref{coromain}}
To control the bias term, we use the following lemma

\begin{lem} \label{approx}Let  $\pi_A$ belong to $B_{2,\infty}^{\boldsymbol{\alpha}}(A)$.
We consider that $S_m'$ is one of the following spaces on $A$
:\begin{itemize}
\item a space of piecewise polynomials of degrees bounded by $s_i>\alpha_i -1$ ($i=1,2$)
based on a partition with rectangles of vertices $1/D_{m_1}$ and $1/D_{m_2}$,
\item a linear span of $\{\phi_\lambda\psi_\mu, \lambda\in \cup_0^{m_1}\Lambda(j), \mu\in \cup_0^{m_2}M(k)\}$
where $\{\phi_\lambda\}$ and $\{\psi_\mu\}$ are orthonormal wavelet bases of respective 
regularities $s_1>\alpha_1 -1$ and $s_2>\alpha_2 -1$ (here $D_{m_i}=2^{m_i}, i=1,2$),
\item the space of trigonometric polynomials with degree smaller than $D_{m_1}$ in the first direction
and smaller than $D_{m_2}$ in the second direction.
\end{itemize}
Let $\pi_m'$ be the orthogonal projection of $\pi_A$ on $S_m'$. 
Then, there exists a positive constant $C_0$ such that
$$\left(\int_A|\pi_A-\pi_m'|^2\right)^{1/2}\leq C_0 [ D_{m_1}^{-\alpha_1} +D_{m_2}^{-\alpha_2} ].$$
\end{lem}

\emph{Proof:}
It is proved in \cite{Hochmuth02} for $S_m'$ a space of wavelets or polynomials and in \cite{Nikolskij}
(p. 191 and 200) for a space of trigonometric polynomials that
$$\left(\int_A|\pi_A-\pi_m'|^2\right)^{1/2}\leq C[\omega_{s_1+1,1}(\pi,D_{m_1}^{-1})+\omega_{s_2+1,2}(\pi,D_{m_2}^{-1})].$$
The definition of $B_{2,\infty}^{\boldsymbol{\alpha}}(A)$ implies
$\left(\int_A|\pi_A-\pi_m'|^2\right)^{1/2}\leq C_0 [D_{m_1}^{-\alpha_1} +D_{m_2}^{-\alpha_2}].$
\hfill$\Box$

If we choose for $S_m'$ the set of the restrictions to $A$ of the functions of $S_m$ and $\pi_A$
the restriction of $\pi$ to $A$, we can apply Lemma \ref{approx}. But $\pi_m'$ is also the restriction to $A$ of $\pi_m$ so that 
$$\|\pi\1_A-\pi_m\|\leq C_0 [ D_{m_1}^{-\alpha_1} +D_{m_2}^{-\alpha_2} ].$$
According to  Theorem \ref{main}
$$\mathbb{E}\|\tilde{\pi}-\pi\1_A\|_n^2 \leq C''\underset{m\in\mathcal{M}_n}{\inf}\left\{D_{m_1}^{-2\alpha_1} 
    +D_{m_2}^{-2\alpha_2}+\frac{D_{m_1}D_{m_2}}{n}\right\}.$$
In particular, if $m^*$ is such that 
$D_{m_1^*}=\lfloor n^\frac{\alpha_2}{\alpha_1+\alpha_2+2\alpha_1\alpha_2}\rfloor$ and
$D_{m_2^*}=\lfloor (D_{m_1^*})^{\frac{\alpha_1}{\alpha_2}}\rfloor$
then $$\mathbb{E}\|\tilde{\pi}-\pi\1_A\|_n^2 \leq 
C'''\left\{D_{m_1^*}^{-2\alpha_1} +\frac{D_{m_1^*}^{1+\alpha_1/\alpha_2}}{n}\right\}
=O\left(n^{-\frac{2\alpha_1\alpha_2}{\alpha_1+\alpha_2+2\alpha_1\alpha_2}}\right).$$
But the harmonic mean of $\alpha_1$ and $\alpha_2$ is 
$\bar\alpha={2\alpha_1\alpha_2}/({\alpha_1+\alpha_2}).$
Then $\mathbb{E}\|\tilde{\pi}-\pi\1_A\|_n^2 =O(n^{-\frac{2\bar\alpha}{2\bar\alpha+2}}).$

The condition $D_{m_1}\leq n^{1/3}$  allows this choice of $m$ only if 
$\frac{\alpha_2}{\alpha_1+\alpha_2+2\alpha_1\alpha_2}<\frac{1}{3}$ i.e. if $\alpha_1-2\alpha_2+2\alpha_1\alpha_2>0$.
In the same manner, the condition $\alpha_2-2\alpha_1+2\alpha_1\alpha_2>0$ must be verified.

\subsection{Proof of Theorem \ref{main2}}

We use the same notations as for the proof of Theorem \ref{main}.
Let us write $$\E\|\tilde{\pi}^*-\pi\1_A\|^2=B_1+B_2+B_3$$ with 
$\begin{cases}
\displaystyle   B_1=\E\left(\|\tilde{\pi}^*-\pi\1_A\|^2\1_{\Omega_\rho^*}\1_{\|\tilde{\pi}\|\leq k_n}\right)\\
\displaystyle B_2=\E\left(\|\tilde{\pi}^*-\pi\1_A\|^2\1_{\Omega_\rho^*}\1_{\|\tilde{\pi}\|> k_n}\right)\\
\displaystyle B_3=\E\left(\|\tilde{\pi}^*-\pi\1_A\|^2\1_{\Omega_\rho^{*c}}\right)\\
\end{cases}$

To bound the first term, we observe that for all $m\in\mathcal{M}_n$, on $\Omega_\rho^*$, 
$\|\tilde{\pi}-\pi_m\|^2\leq f_0^{-1}\rho\|\tilde{\pi}-\pi_m\|_n^2$. Then 
\begin{eqnarray*}
&&  \|\tilde{\pi}-\pi\1_A\|^2\1_{\Omega_\rho^*}\leq 2\|\tilde{\pi}-\pi_m\|^2\1_{\Omega_\rho^*}+
2\|{\pi}_m-\pi\1_A\|^2\\
&&\leq 2f_0^{-1}\rho\|\tilde{\pi}-\pi_m\|_n^2\1_{\Omega_\rho^*}+2\|{\pi}_m-\pi\1_A\|^2\\
&&\leq 2f_0^{-1}\rho\{2\|\tilde{\pi}-\pi\1_A\|_n^2\1_{\Omega_\rho^*}+2\|{\pi}_m-\pi\1_A\|_n^2\}+2\|{\pi}_m-\pi\1_A\|^2
\end{eqnarray*}
Thus
$$B_1\leq\E\left(\|\tilde{\pi}-\pi\1_A\|^2\1_{\Omega_\rho^*}\right)\leq
4f_0^{-1}\rho\E\left(\|\tilde{\pi}-\pi\1_A\|_n^2\1_{\Omega_\rho^*}\right)
+(4f_0^{-1}\rho\|f\|_\infty+2)\|{\pi}_m-\pi\1_A\|^2.$$
But, using \eqref{firstterm}, we obtain
$$B_1\leq (24f_0^{-1}\rho\|f\|_\infty+2)
\|\pi_m-\pi\1_{A}\|^2+24f_0^{-1}\rho\pen(m)+36f_0^{-1}\rho^2\frac{C_1}{n}.$$
Since $\rho=3/2$, by setting $C=\max(36f_0^{-1}\|f\|_\infty+1,36f_0^{-1})$,
$$B_1\leq C\{\|\pi_m-\pi\1_{A}\|^2+\pen(m)\}+\frac{81f_0^{-1}C_1}{n}$$
for all $m\in\mathcal{M}_n$.

Next, the definition of $\tilde{\pi}^*$ and the Markov inequality provide
\begin{equation}\label{B2}
B_2\leq\E\left(\|\pi\1_A\|^2\1_{\Omega_\rho^*}\1_{\|\tilde{\pi}\|> k_n}\right)
\leq \|\pi\|^2\frac{\E(\|\tilde{\pi}\|^2\1_{\Omega_\rho^*})}{k_n^2}.
\end{equation}
But $\|\tilde{\pi}\|^2\1_{\Omega_\rho^*}\leq \rho f_0^{-1} \|\tilde{\pi}\|_n^2
\leq 2\rho f_0^{-1} (\|\tilde{\pi}-\pi\1_A\|_n^2+\|\pi\1_A\|_n^2).$
Now we use \eqref{majdepi-pichapeau} to state
\begin{eqnarray*}
   \|\tilde{\pi}\|^2\1_{\Omega_\rho^*}&\leq& 2\rho f_0^{-1} (\|\pi\|_\infty+4\phi_2n^{1/3}+\|\pi\1_A\|_n^2)\\
&\leq& 2\rho f_0^{-1} (\|\pi\|_\infty+4\phi_2n^{1/3}+\frac{1}{n}\sum_{i=1}^n\|\pi\|_\infty\int\pi(X_i,y)dy]\\
&\leq& 2\rho f_0^{-1} (2\|\pi\|_\infty+4\phi_2n^{1/3}).
\end{eqnarray*}
Then, since $k_n=n^{2/3}$, \eqref{B2} becomes
$$B_2\leq\|\pi\|^2\frac{2\rho f_0^{-1} (2\|\pi\|_\infty+4\phi_2n^{1/3})}{k_n^2}
\leq 4\rho f_0^{-1}\|\pi\|^2\left(\frac{\|\pi\|_\infty}{n^{4/3}}+\frac{2\phi_2}{n}\right).$$

Lastly 
$$B_3\leq \E\left(2(\|\tilde{\pi}^*\|^2+\|\pi\1_A\|^2)\1_{\Omega_\rho^{*c}}\right)
\leq 2(k_n^2+\|\pi\|^2)P(\Omega_\rho^{*c}).$$
We now remark that 
$P(\Omega_\rho^{*c})=P(\Omega^{*c})+P(\Omega_\rho^c\cap\Omega^*).$
In the geometric case $\beta_{2q_n}\leq e^{-\gamma c\ln(n)}\leq n^{-\gamma c}$ and in the other case
$\beta_{2q_n}\leq (2q_n)^{-\gamma}\leq n^{-\gamma c}$. 
Then $$P(\Omega^{*c})\leq 4p_n\beta_{2q_n}\leq n^{1-c\gamma}.$$ 
But, if $\gamma>20$ in the arithmetic case, we can choose $c$ such that $c\gamma>\cfrac{10}{3}$ 
and so $P(\Omega^{*c})\leq n^{-7/3}.$
Then, using Proposition \ref{omegac},
$$B_3\leq2(n^{4/3}+\|\pi\|^2)\frac{1+C_2}{n^{7/3}}\leq \frac{2(C_2+1)(1+\|\pi\|^2)}{n}.$$

\subsection{Proof of Theorem \ref{lowerbound}}
Let $\psi$ be a very regular wavelet with compact support.
For $J=(j_1,j_2)\in \mathbb{Z}^2$ to be chosen below and $K=(k_1, k_2)\in \mathbb{Z}^2$, we set
$$\psi_{JK}(x,y)=2^{(j_1+j_2)/2}\psi(2^{j_1}x-k_1)\psi(2^{j_2}y-k_2).$$
Let $\pi_0(x,y)=c_0\1_{B}(y)$ with $B$ a compact set such that
$A\subset B\times B$ and $|B|\geq 2|A|^{1/2}/L$, and $c_0=|B|^{-1}$.
So $\pi_0$ is a transition density with $\|\pi_0\|_{B_{2,\infty}^{\boldsymbol{\alpha}}(A)}\leq L/2.$
Now we set $R_J$ the maximal subset of $\mathbb{Z}^2$ such that 
\begin{equation*}
  \text{Supp} (\psi_{JK})\subset A\quad   \forall K\in R_J, \qquad
  \text{Supp} (\psi_{JK}) \cap \text{Supp} (\psi_{JK'})=\emptyset \text{ if } K\neq K'.
\end{equation*}
The cardinal of $R_J$ is $|R_J|=c2^{j_1+j_2}$, with $c$ a positive constant which depends only on $A$
and the support of $\psi$.
Let, for all $\varepsilon=(\varepsilon_K)\in \{-1,1\}^{|R_J|}$,
$$\pi_\varepsilon=\pi_0+\frac{1}{\sqrt{n}}\sum_{K\in R_J}\varepsilon_K\psi_{JK}.$$
Let us denote by $\mathcal{G}$ the set of all such $\pi_\varepsilon$.
Since $\int \psi=0$ and $\pi_0$ is a transition density, for all $x$ in $\mathbb{R}$,
$\int\pi_\varepsilon(x,y)dy=1.$
Additionally $\pi_\varepsilon(x,y)=\pi_0(x,y)\geq 0$ if $(x,y)\notin A$, and if $(x,y)\in A$:
$\pi_\varepsilon \geq c_0- 2^{(j_1+j_2)/2}\|\psi\|_\infty^2/{\sqrt{n}}$
and then  $\pi_\varepsilon(x,y)\geq c_0/2>0$ as soon as
\begin{equation}
  \left(\frac{2^{j_1+j_2}}{n}\right)^{1/2}\leq \frac{c_0}{2\|\psi\|_\infty^2}.
\label{condition1}\end{equation}
Thus, if \eqref{condition1} holds, $\pi_\varepsilon(x,y) \geq (c_0/2)\1_B(y)$ for all $x, y$.
It implies that the underlying Markov chain is Doeblin recurrent 
and then positive recurrent. We verify that $f=c_0\1_{B}$ is the stationary density.
To prove that $\pi_\varepsilon\in \mathcal{B}$, we still have to compute $\|\pi_\varepsilon\|_{B_{2,\infty}^{\boldsymbol{\alpha}}(A)}$.
Hochmuth \cite{Hochmuth02} proves that for $\psi$ smooth enough 
$\|\sum_{K\in R_J}\varepsilon_K\psi_{JK}\|_{B_{2,\infty}^{\boldsymbol{\alpha}}(A)} \leq 
(2^{j_1\alpha_1}+2^{j_2\alpha_2})\|\sum_{K\in R_J}\varepsilon_K\psi_{JK}\|_{A}$. Since
$$\|\sum_{K\in R_J}\varepsilon_K\psi_{JK}\|_{A}^2  =\sum_{K\in R_J}|\varepsilon_K|^2
= c2^{j_1+j_2},$$
  then 
$$\|\pi_\varepsilon\|_{B_{2,q}^{\boldsymbol{\alpha}}(A)}\leq \frac{L}{2}+
\frac{2^{j_1\alpha_1}+2^{j_2\alpha_2}}{\sqrt{n}}c^{1/2}2^{(j_1+j_2)/2}.$$
From now on, we suppose that Condition C is verified where 
$$\text{Condition C:  }\qquad
\frac{(2^{j_1\alpha_1}+2^{j_2\alpha_2})2^{(j_1+j_2)/2}}{\sqrt{n}}\leq \frac{L}{2c^{1/2}}.$$
It entails in particular that \eqref{condition1} holds if $j_1$ and $j_2$ are great enough.
Then for all $\varepsilon$, $\pi_\varepsilon\in \mathcal{B}$.
We now use the Lemma 10.2 p.160 in H\"{a}rdle et al. \cite{HKPT}. 
The likelihood ratio can be written 
$$\Lambda_n(\pi_{\varepsilon_{*K}},\pi_\varepsilon)=\prod_{i=1}^{n}\frac{\pi_{\varepsilon_{*K}}(X_i,X_{i+1})}
{\pi_\varepsilon(X_i,X_{i+1})}.$$ Note that $\pi_\varepsilon(X_i,X_{i+1})>0$ $P_{\pi_\varepsilon}$- and 
$P_{\pi_{\varepsilon_{*K}}}$- almost surely (actually the chain ``lives'' on $B$). Then 
$$\log(\Lambda_n(\pi_{\varepsilon_{*K}},\pi_\varepsilon))=\sum_{i=1}^{n}\log\left(1-\frac{2}{\sqrt{n}}
\frac{\varepsilon_K\psi_{JK}(X_i,X_{i+1})}{\pi_\varepsilon(X_i,X_{i+1})}\right)$$
We set $U_{JK}(X_i,X_{i+1})=-\varepsilon_K\psi_{JK}(X_i,X_{i+1})/{\pi_\varepsilon(X_i,X_{i+1})}$ so that

\begin{eqnarray*}
&&  \log(\Lambda_n(\pi_{\varepsilon_{*K}},\pi_\varepsilon))=\sum_{i=1}^{n}\log
\left(1+\frac{2}{\sqrt{n}}U_{JK}(X_i,X_{i+1})\right)\\
&&=\sum_{i=1}^{n}\left\{\theta\left(\frac{2}{\sqrt{n}}U_{JK}(X_i,X_{i+1})\right)
+\frac{2}{\sqrt{n}}U_{JK}(X_i,X_{i+1})-\frac{2}{n}U_{JK}^2(X_i,X_{i+1})\right\}\\
&&=u_n+v_n-w_n
\end{eqnarray*}
with $\theta$ the function defined by $\theta(u)=\log(1+u)-u+\cfrac{u^2}{2}$.
Now we prove the three following assertions
\begin{description}
\item[$1^{\circ}$] $\E_{\pi_\varepsilon}(|u_n|)=\E_{\pi_\varepsilon}\left(\left|\sum_{i=1}^{n}
\theta\left(\frac{2}{\sqrt{n}}U_{JK}(X_i,X_{i+1})\right)\right|\right)\underset{n\to\infty}{\to} 0$
\item[$2^{\circ}$] $\E_{\pi_\varepsilon}(w_n)=\E_{\pi_\varepsilon}\left(\cfrac{2}{n}\sum_{i=1}^{n}U_{JK}^2(X_i,X_{i+1})\right)\leq 4$
\item[$3^{\circ}$] $\E_{\pi_\varepsilon}(v_n^2)=\E_{\pi_\varepsilon}\left(\cfrac{4}{n}|\sum_{i=1}^{n}U_{JK}(X_i,X_{i+1})|^2\right) \leq 8$
\end{description}

$1^{\circ}:$ First we observe that $\left\|\cfrac{2}{\sqrt{n}}U_{JK}\right\|_\infty\leq \cfrac{2}{\sqrt{n}}
\cfrac{2^{(j_1+j_2)/2}\|\psi\|_\infty^2}{c_0/2}=O\left(\cfrac{2^{(j_1+j_2)/2}}{\sqrt{n}}\right)$
and $\cfrac{2^{(j_1+j_2)}}{n}\to 0$ since Condition C holds.
So there exists some integer $n_0$ such that $\forall n\geq n_0$, 
$\forall x, y$, $|\theta(2U_{JK}(x,y)/\sqrt{n})|\leq |2U_{JK}(x,y)/\sqrt{n}|^3.$
But
\begin{eqnarray*}
 \iint\left|\cfrac{2U_{JK}(x,y)}{\sqrt{n}}\right|^3f(x)\pi_\varepsilon(x,y)dxdy
=\cfrac{8}{n\sqrt{n}}\iint\cfrac{|\psi_{JK}(x,y)|^3}{\pi_\varepsilon(x,y)^2}f(x)dxdy\\
\leq \cfrac{8}{n\sqrt{n}}\cfrac{2^{(j_1+j_2)/2}\|\psi\|_\infty^2c_0}{(c_0/2)^2}\iint\psi_{JK}(x,y)^2dxdy
\leq\frac{32\|\psi\|_\infty^2}{c_0n}\left(\cfrac{2^{(j_1+j_2)}}{n}\right)^{1/2}.
\end{eqnarray*}
Then $\displaystyle \E_{\pi_\varepsilon}|u_n|\leq \sum_{i=1}^{n}\cfrac{32\|\psi\|_\infty^2}{c_0n}
\left(\cfrac{2^{(j_1+j_2)}}{n}\right)^{1/2}\to_{n\to\infty} 0.$

$2^{\circ}:$ We bound the expectation of $U_{JK}(X_i,X_{i+1})^2 $: 
\begin{equation}
  \label{aa}
\E_{\pi_\varepsilon}(U_{JK}(X_i,X_{i+1})^2)=\iint \cfrac{\psi_{JK}^2(x,y)}{\pi_\varepsilon(x,y)}f(x)dxdy
\leq c_0\iint_A \cfrac{\psi_{JK}^2(x,y)}{c_0/2}dxdy\leq 2.
\end{equation}
And then $\E_{\pi_\varepsilon}(w_n)=\E_{\pi_\varepsilon}\left((2/n)\sum_{i=1}^{n}U_{JK}(X_i,X_{i+1})^2\right)
\leq 4$.

$3^{\circ}:$ We observe that $\E_{\pi_\varepsilon}(U_{JK}(X_i,X_{i+1})|X_1,\dots,X_i)= 0$
and thus $\sum_{i=1}^nU_{JK}(X_i,X_{i+1})$ is a martingale.
A classic property of square integrable martingales involves
$$E_{\pi_\varepsilon}\left(\left[\sum_{i=1}^n U_{JK}(X_i,X_{i+1})\right]^2\right)=
\sum_{i=1}^n\E_{\pi_\varepsilon}\left[U_{JK}(X_i,X_{i+1})^{2}\right].$$
\label{mart}
Thus, using \eqref{aa},
$\E_{\pi_\varepsilon}(v_N^2)=(4/n)\sum_{i=1}^n\E_{\pi_\varepsilon}\left[U_{JK}(X_i,X_{i+1})^{2}\right]\leq 8.$

We deduce easily from the three previous assertions $1^{\circ}$, $2^{\circ}$ and $3^{\circ}$ that 
there exists $\lambda>0$ and $p_0$ such that
$P_{\pi_\varepsilon}(\Lambda_n(\pi_{\varepsilon_{*K}},\pi_\varepsilon)>e^{-\lambda})\geq p_0$.
Thus, 
according to Lemma 10.2 in \cite{HKPT},
$$\max_{\pi_\varepsilon\in\mathcal{G}}\E_{\pi_\varepsilon}\|\hat{\pi}_n-\pi_\varepsilon\|_{A}^2\geq \frac{|R_J|}{2}
\delta^2e^{-\lambda}p_0$$
where 
$\delta = \inf_{\varepsilon\neq\varepsilon'}\|\pi_\varepsilon-\pi_{\varepsilon'}\|_{A}/2=
\left\|\varepsilon_K\psi_{JK}/{\sqrt{n}}\right\|_{A}={1}/{\sqrt{n}}.$

Now for all $n$ we choose $J=J(n)=(j_1(n),j_2(n))$ such that
$$c_1/2\leq 2^{j_1}n^{-\frac{\alpha_2}{\alpha_1+\alpha_2+2\alpha_1\alpha_2}}\leq c_1
\quad \text{ and }\quad
c_2/2\leq 2^{j_2}n^{-\frac{\alpha_1}{\alpha_1+\alpha_2+2\alpha_1\alpha_2}}\leq c_2$$
with $c_1$ and $c_2$ such that  
$(c_1^{\alpha_1}+c_2^{\alpha_2})\sqrt{c_1c_2}\leq {L}/(2c^{1/2})$
so that Condition C is satisfied. 
Moreover, we have
$$|R_J|\delta^2\geq\frac{cc_1c_2}{4}n^{\frac{\alpha_2+\alpha_1}{\alpha_1+\alpha_2+2\alpha_1\alpha_2}-1}
\geq\frac{cc_1c_2}{4}n^{\frac{-2\alpha_1\alpha_2}{\alpha_1+\alpha_2+2\alpha_1\alpha_2}}$$
Thus 
$$\max_{\pi_\varepsilon\in\mathcal{G}}\E_{\pi_\varepsilon}\|\hat{\pi}_n-\pi_\varepsilon\|_{A}^2\geq
\frac{ce^{-\lambda}p_0c_1c_2}{8}n^{\frac{-2\alpha_1\alpha_2}{\alpha_1+\alpha_2+2\alpha_1\alpha_2}}.$$
And then for all estimator
$$\sup_{\pi \in \mathcal{B}}\E_\pi\|\hat{\pi}_n-\pi\|_{A}^2\geq Cn^{-\frac{2\bar\alpha}{2\bar\alpha+2}}$$
with $C=ce^{-\lambda}p_0c_1c_2/8$.

\subsection{Proof of Proposition \ref{talagrand}}

Let $\begin{cases}
\Gamma_i(t)=t(X_i,X_{i+1})-\int t(X_i,y)\pi(X_i,y)dy,\\
\Gamma_i^*(t)=t(X_i^*,X_{i+1}^*)-\int t(X_i^*,y)\pi(X_i^*,y)dy,\\
\Gamma_i^{**}(t)=t(X_i^{**},X_{i+1}^{**})-\int t(X_i^{**},y)\pi(X_i^{**},y)dy.\\
\end{cases}$

We now define $Z_n^*(t)$:
$$\displaystyle Z_n^*(t)=\frac{1}{n}\sum_{i \text{ odd}}\Gamma_i^*(t)+
\frac{1}{n}\sum_{i\text{ even}}\Gamma_i^{**}(t).$$
Let us remark that $Z_n^*(t)\1_{\Omega^*}=Z_n(t)\1_{\Omega^*}.$
Next we split each of these terms :

$$\displaystyle Z_{n,1}^*(t)=\frac{1}{n}\sum_{l=0}^{p_n-1}
\sum_{i=4lq_n+1, i \text{ odd}}^{2(2l+1)q_n-1}\Gamma_i^*(t), \quad
\displaystyle Z_{n,2}^*(t)=\frac{1}{n}\sum_{l=0}^{p_n-1}
\sum_{i=2(2l+1)q_n+1, i \text{ odd}}^{2(2l+2)q_n-1}\Gamma_i^*(t),$$
$$\displaystyle Z_{n,3}^*(t)=\frac{1}{n}\sum_{l=0}^{p_n-1}
\sum_{i=4lq_n+2, i \text{ even}}^{2(2l+1)q_n}\Gamma_i^{**}(t), \quad
\displaystyle Z_{n,4}^*(t)=\frac{1}{n}\sum_{l=0}^{p_n-1}
\sum_{i=2(2l+1)q_n+2, i \text{ even}}^{2(2l+2)q_n}\Gamma_i^{**}(t).$$

We use the following lemma:

\begin{lem}(Talagrand \cite{talagrand1996})\\
Let $\mathcal{U}_0,\dots,\mathcal{U}_{N-1}$ i.i.d. variables and $(\zeta_t)_{t \in B}$ a set of functions. \\
Let $\displaystyle G(t)= \cfrac{1}{N}\sum_{l=0}^{N-1} \zeta_t(\mathcal{U}_l).$
We suppose that \\
(1) $\underset{t\in B}{\sup}\|\zeta_t\|_\infty \leq M_1, \quad$
(2) $\mathbb{E}(\underset{t\in B}{\sup}|G(t)|)\leq H, \quad$
(3) $\underset{t\in B}{\sup}Var [\zeta_t(\mathcal{U}_0)]\leq v$.

Then, there exists $K>0$, $K_1>0$, $K_2>0$ such that 
$$\mathbb{E}\left[\underset{t\in B}{\sup}G^2(t)-10H^2\right]_+
\leq K\left[\frac{v}{N}e^{-K_1\frac{NH^2}{v}}+
\frac{M_1^2}{N^2}e^{-K_2\frac{NH}{M_1}}\right]$$
\label{tala}\end{lem}

Here $N=p_n$, $B=B_f(m')$ and for $l \in \{0,\dots,p_n-1\}$,
$\mathcal{U}_l=(X^*_{4lq_n+1},..,X^*_{2(2l+1)q_n}),$
$$\zeta_t(x_1,\dots,x_{2q_n})=\frac{1}{q_n}\sum_{i=1, i \text{ odd }}^{2q_n-1}
t(x_i,x_{i+1})-\int t(x_i,y)\pi(x_i,y)dy.$$
Then 
$$ G(t)= \cfrac{1}{p_n}\sum_{l=0}^{p_n-1} \frac{1}{q_n}\sum_{i=4lq_n+1, i\text{ odd }}^{2(2l+1)q_n-1}
\Gamma_i^*(t)=4Z_{n,1}^*(t).$$
We now compute $M_1$, $H$ and $v$.
 
(1)We recall that $S_m+S_{m'}$ is included in the model $S_{m''}$ 
with dimension $\max(D_{m_1}, D_{m_1'})\max(D_{m_2}, D_{m_2'})$. 
\begin{eqnarray*}
\underset{t\in B}{\sup}\|\zeta_t\|_\infty &\leq &\underset{t\in B}{\sup}\|t\|_\infty
\frac{1}{q_n}\sum_{i=1, i \text{ odd }}^{2q_n-1}\left(1+\int \pi(x_i,y)dy\right)\\
&\leq &2\phi_0\sqrt{\max(D_{m_1}, D_{m_1'})\max(D_{m_2}, D_{m_2'})}\|t\|\leq \frac{2\phi_0}{f_0} n^{1/3}.
\end{eqnarray*}
Then we set $M_1=\cfrac{2\phi_0}{f_0}n^{1/3}  $.

(2) Since $A_0 \text{ and } A_0^*$ have the same distribution, $\zeta_t(\mathcal{U}_0)=\frac{1}{q_n}
\sum_{i=1, i\text{ odd }}^{2q_n-1}\Gamma_i^*(t)$ has the same distribution than 
$\frac{1}{q_n}\sum_{i=1, i\text{ odd }}^{2q_n-1}\Gamma_i(t)$.
We observe that $\E(\Gamma_i(t)|X_i)=0$ and then for all set $I$
\begin{eqnarray*}
&&\E\left(\left[\sum_{i\in I} \Gamma_i(t)\right]^2\right)=\E\left(\sum_{i,j\in I}\Gamma_i(t)\Gamma_j(t)\right)\\
&&=2\E\left(\sum_{j<i}\E[\Gamma_i(t)\Gamma_j(t)|X_i]\right)+\sum_{i\in I}\E\left[\Gamma_i^{2}(t)\right]\\
&&=2\E\left(\sum_{j<i}\Gamma_j(t)\E[\Gamma_i(t)|X_i]\right)+\sum_{i\in I}\E\left[\Gamma_i^{2}(t)\right]
=\sum_{i\in I}\E\left[\Gamma_i^{2}(t)\right].
\end{eqnarray*}
In particular 
\begin{eqnarray*}
{\rm Var} [\zeta_t(\mathcal{U}_0)]&=&\E\left(\left[\frac{1}{q_n}\sum_{i=1, i\text{ odd }}^{2q_n-1}
\Gamma_i(t)\right]^2\right)
=\frac{1}{q_n^2}\sum_{i=1, i\text{ odd }}^{2q_n-1}\E\left[\Gamma_i^{2}(t)\right]\\
&\leq&\frac{1}{q_n^2}\sum_{i=1, i\text{ odd }}^{2q_n-1}\E\left[t^2(X_i,X_{i+1})\right]
\leq\frac{1}{q_n}\|\pi\|_\infty\|t\|_f^2.
\end{eqnarray*}
Then $v=\cfrac{\|\pi\|_\infty}{q_n}$.
 
(3) Let $(\bar\varphi_j\otimes\psi_k)_{(j,k)\in \Lambda(m,m')}$ an orthonormal basis of 
$(S_m +S_{m'}, \|.\|_f).$ 
\begin{eqnarray*}
\E(\underset{t\in B}{\sup}|G^2(t)|)&\leq& \sum_{j,k}\E(G^2(\bar\varphi_j\otimes\psi_k))\\
&\leq&\sum_{j,k}\cfrac{1}{p_n^2q_n^2}\E\left(\left[\sum_{l=0}^{p_n-1} 
\sum_{i=4lq_n+1, i\text{ odd }}^{2(2l+1)q_n-1}\Gamma_i^*(\bar\varphi_j\otimes\psi_k)\right]^2\right)\\
&\leq&\sum_{j,k}\cfrac{16}{n^2}\sum_{l=0}^{p_n-1}\E\left(\left[\sum_{i=4lq_n+1, i\text{ odd }}^{2(2l+1)q_n-1}
\Gamma_i^*(\bar\varphi_j\otimes\psi_k)\right]^2\right)
\end{eqnarray*}
where we used the independence of the $A_l^*$. Now we can replace $\Gamma_i^*$ by $\Gamma_i$ in 
the sum because $A_l$ and $A_l^*$ have the same distribution and we use as previously the martingale 
property of the $\Gamma_i$.

\begin{eqnarray*}
\E(\underset{t\in B}{\sup}|G^2(t)|)&\leq&\sum_{j,k}\cfrac{16}{n^2}
\sum_{l=0}^{p_n-1}\E\left(\left[\sum_{i=4lq_n+1, i\text{ odd }}^{2(2l+1)q_n-1}
\Gamma_i(\bar\varphi_j\otimes\psi_k)\right]^2\right)\\
&\leq&\sum_{j,k}\cfrac{16}{n^2}\sum_{l=0}^{p_n-1}\sum_{i=4lq_n+1, i\text{ odd }}^{2(2l+1)q_n-1}
\E\left(\Gamma_i^2(\bar\varphi_j\otimes\psi_k)\right)\\
&\leq&\sum_{j,k}\cfrac{4}{n}\|\pi\|_\infty\|\bar\varphi_j\otimes\psi_k\|_f^2
\leq 4\|\pi\|_\infty\cfrac{D(m,m')}{n}.
\end{eqnarray*}
Then $\E^2(\underset{t\in B}{\sup}|G(t)|)\leq4\|\pi\|_\infty\cfrac{D(m,m')}{n}$ and 
$H^2=4\|\pi\|_\infty\cfrac{D(m,m')}{{n}}.$

According to Lemma  \ref{tala}, there exists $K'>0$, $K_1>0$, $K'_2>0$ such that 
 $$\E\left[\underset{t\in B_f(m')}{\sup}(4Z_{n,1}^*)^2(t)-10H^2\right]_+
 \leq K'\left[\frac{1}{n}e^{-K_1D(m,m')}+n^{-4/3}q_n^2e^{-K'_2n^{1/6}\sqrt{D(m,m')}/q_n}\right].$$

But $q_n\leq n^c$ with $c<\frac{1}{6}$. So
\begin{eqnarray}\label{eqta}
&& \sum_{m'\in \mathcal{M}_n}\E\left[\underset{t\in B_f(m')}{\sup}Z_{n,1}^{*2}(t)-\frac{p(m,m')}{4}\right]_+\nonumber\\
&& \leq \frac{K'}{n}\left[\sum_{m'\in \mathcal{M}_n}e^{-K_1D(m,m')}+
n^{2c-1/3}|\mathcal{M}_n|e^{-K'_2n^{1/6-c}}\right]\leq \frac{A_1}{n}.
\end{eqnarray}

In the same way,
$\sum_{m'\in \mathcal{M}_n}\E\left[\underset{t\in B_f(m')}{\sup}Z_{n,r}^{*2}(t)-{p(m,m')}/{4}\right]_+
\leq {A_r}/{n}$ for $r=2,3,4.$
And then 
$$\begin{array}{l}
\displaystyle\sum_{m'\in \mathcal{M}_n}\E\left(\left[\underset{t\in B_f(m')}{\sup}Z_{n}^{2}(t)-p(m,m')\right]_+
\1_{\Omega^*}\right)\\\qquad
\displaystyle =\sum_{m'\in \mathcal{M}_n}\E\left(\left[\underset{t\in B_f(m')}{\sup}Z_{n}^{*2}(t)-p(m,m')\right]_+
\1_{\Omega^*}\right)\leq \frac{C_1}{n}.
\end{array}$$

\subsection{Proof of Proposition \ref{omegac}}

First we observe that
$$\displaystyle P(\Omega_\rho^c\cap\Omega^*)\leq 
P\left(\sup_{t\in \mathcal{B}} |\nu_n(t^2)|>1-1/\rho\right)$$
where $\displaystyle\nu_n(t)=\frac{1}{n}\sum_{i=1}^{n}\int[t(X_i^*,y)-\E(t(X_i^*,y))]dy$
and $\mathcal{B}=\{t \in \mathcal{S} \quad \|t\|_f=1\}.$

\noindent But, if $t(x,y)=\sum_{j,k}a_{j,k}\varphi_j(x)\psi_k(y)$, then
$$ \nu_n(t^2)= \sum_{j,j'}\sum_{k}a_{j,k}a_{j',k} \bar\nu_n(\varphi_j\varphi_{j'})$$
where \begin{equation}
\displaystyle\bar\nu_n(u)=\frac{1}{n}\sum_{i=1}^{n}[u(X_i^*)-\E(u(X_i^*))].
\label{nubar}\end{equation} 
Let $b_j=(\sum_ka_{j,k}^2)^{1/2}$, then $ |\nu_n(t^2)|\leq \sum_{j,j'}b_jb_{j'}|\bar\nu_n(\varphi_j\varphi_{j'})|$
and, if $t \in \mathcal{B}$, ${\sum_j b_j^2=\sum_j\sum_k a_{j,k}^2=\|t\|^2\leq f_0^{-1}}$.

Thus
$$\sup_{t\in \mathcal{B}} |\nu_n(t^2)|\leq f_0^{-1}\sup_{\sum b_{j}^2=1}\sum_{j,l}b_{j}b_{l}
|\bar\nu_n(\varphi_j\varphi_{l})|.$$

\begin{lem}
Let $B_{j,l}=\|\varphi_j\varphi_{l}\|_\infty$ and $V_{j,l}=\|\varphi_j\varphi_{l}\|_2$.
Let, for any symmetric matrix $(A_{j,l})$
$$ \bar\rho(A)=\sup_{\sum a_j^2=1} \sum_{j,l}|a_ja_{l}|A_{j,l}$$
and $L(\varphi)=\max\{\bar\rho^2(V),\bar\rho(B)\}.$
Then, if M2 is satisfied, $ L(\varphi)\leq\phi_1\mathcal{D}_n^2$.
\end{lem}

This lemma is proved in Baraud et al. \cite{BCV}.

Let $x=\cfrac{f_0^2(1-1/\rho)^2}{40\|f\|_\infty L(\varphi)}$ and 
$\Delta=\left\{\forall j \forall l \quad|\bar\nu_n(\varphi_j\varphi_l)|\leq
        4\left[B_{j,l} x+V_{j,l}\sqrt{2\|f\|_\infty x}\right]\right\}. $
On $\Delta$:
\begin{eqnarray*}
\sup_{t\in \mathcal{B}} |\nu_n(t^2)|
&\leq&  4f_0^{-1}\sup_{\sum b_{j}^2=1}\sum_{j,l}b_{j}b_{l}\left[B_{j,l} x+V_{j,l}\sqrt{2\|f\|_\infty x}\right]\\
&\leq&  4f_0^{-1}\left[\bar\rho(B) x+\bar\rho(V)\sqrt{2\|f\|_\infty x}\right]\\
&\leq&(1-1/\rho)\left[\frac{f_0(1-1/\rho)}{10\|f\|_\infty}\frac{\bar\rho(B)}{L(\varphi)}
+\frac{2}{\sqrt{5}}\left(\frac{\bar\rho^2(V)}{L(\varphi)}\right)^{1/2}\right]\\
&\leq&(1-1/\rho)\left[\frac{1}{10}+\frac{2}{\sqrt{5}}\right]\leq(1-1/\rho).
\end{eqnarray*}

Then  $\displaystyle P\left(\sup_{t\in \mathcal{B}} |\nu_n(t^2)|>1-\frac{1}{\rho}\right)\leq P(\Delta^c).$
But $\bar\nu_n(u)=2\bar\nu_{n,1}(u)+2\bar\nu_{n,2}(u)$ with 
$$ \bar\nu_{n,r}(u)=\frac{1}{p_n}\sum_{l=0}^{p_n-1}Y_{l,r}(u)\qquad r=1,2$$
with $\begin{cases}
Y_{l,1}(u)&=\cfrac{1}{2q_n}\sum_{i=4lq_n+1}^{2(2l+1)q_n}[u(X_i^*)-\E(u(X_i^*))],\\
Y_{l,2}(u)&=\cfrac{1}{2q_n}\sum_{i=2(2l+1)q_n+1}^{2(2l+2)q_n}[u(X_i^*)-\E(u(X_i^*))].
\end{cases}$

To bound $P(\bar\nu_{n,1}(\varphi_j\varphi_l)\geq B_{j,l} x+V_{j,l}\sqrt{2\|f\|_\infty x})$, we will use 
the Bernstein inequality given in Birg\'e and Massart \cite{birge&massart98}.
That is why we bound $\E |Y_{l,1}(u)|^m$:
\begin{eqnarray*}
\E |Y_{l,1}(u)|^m&\leq& \cfrac{1}{4q_n^2}(2\|u\|_\infty)^{m-2}
\E\left|\sum_{i=4lq_n+1}^{2(2l+1)q_n}[u(X_i^*)-\E(u(X_i^*))]\right|^2\\
&\leq&(2\|u\|_\infty)^{m-2}\cfrac{1}{4q_n^2}
\E\left|\sum_{i=4lq_n+1}^{2(2l+1)q_n}[u(X_i)-\E(u(X_i))]\right|^2\\
&\leq& (2\|u\|_\infty)^{m-2}\cfrac{1}{4q_n^2}
\E\left|\sum_{i=2lq_n+1}^{2(2l+1)q_n}[u(X_1)-\E(u(X_1))]\right|^2
\end{eqnarray*}
since $ X_i^*=X_i$ on $\Omega^*$ and the $X_i$ have the same distribution than $X_1$. Thus
\begin{eqnarray}\hspace{-0.7cm}
\E |Y_{l,1}(u)|^m&\leq& (2\|u\|_\infty)^{m-2}\E|u(X_1)-\E(u(X_1))|^2
\leq (2\|u\|_\infty)^{m-2}\int u^2(x)f(x)dx\nonumber\\
&\leq &2^{m-2}(\|u\|_\infty)^{m-2}(\sqrt{\|f\|_\infty}\|u\|)^2.\label{bernstein2}
\end{eqnarray}
With $u=\varphi_j\varphi_{j'}$,
$ \E |Y_{l,1}(\varphi_j\varphi_{j'})|^m\leq 2^{m-2}(B_{j,j'})^{m-2}(\sqrt{\|f\|_\infty}V_{j,j'})^2.$
And then
$$P(|\bar\nu_{n,r}(\varphi_j\varphi_l)|\geq B_{j,l} x+V_{j,l}\sqrt{2\|f\|_\infty x})
\leq 2e^{-p_nx}.$$
Given that $P(\Omega_\rho^c\cap\Omega^*)\leq P(\Delta^c)=
\sum_{j,l}P\left(|\bar\nu_n(\varphi_j\varphi_l)|>4(B_{j,l} x+V_{j,l}\sqrt{2\|f\|_\infty x})\right)$,
\begin{eqnarray*}
 P(\Omega_\rho^c\cap\Omega^*)&\leq&4\mathcal{D}_n^2\exp\left\{-\cfrac{p_nf_0^2(1-1/\rho)^2}
{40\|f\|_\infty L(\varphi)}\right\}\\
&\leq&4n^{2/3}\exp\left\{-\cfrac{f_0^2(1-1/\rho)^2}{160\|f\|_\infty} \cfrac{n}{q_nL(\varphi)}\right\}.
\end{eqnarray*}

But $L(\varphi)\leq\phi_1\mathcal{D}_n^2\leq\phi_1 n^{2/3}$ and $q_n\leq n^{1/6}$ so 
\begin{equation}
  \label{eqom}
  P(\Omega_\rho^c\cap\Omega^*)\leq4n^{2/3}\exp\left\{-\cfrac{f_0^2(1-1/\rho)^2}{160\|f\|_\infty\phi_1} n^{1/6}\right\}
\leq \frac{C}{n^{7/3}}.
\end{equation}

\section*{Acknowledgement}
I would like to thank Fabienne Comte for her constructive comments and suggestions, and Marc Hoffmann 
for his help regarding the lower bound.


\section*{Appendix : random penalty}

Here we prove that Theorem \ref{main} is valid with a penalty which does not depend 
on $\|\pi\|_\infty$.
\begin{thm}
We consider the following penalty :
$$\overline{\pen}(m)=\overline{K_0}\|\hat\pi\|_\infty\frac{D_{m_1}D_{m_2}}{n}$$ 
where $\overline{K_0}$ is a numerical constant and 
$\hat\pi=\hat\pi_{m*}$ with $S_{m*}$ a space of trigonometric polynomials such that 
$$\ln n \leq D_{m_1*}=D_{m_2*} \leq n^{1/6}.$$ 
If the restriction of $\pi$ to $A$ belongs to $B_{2,\infty}^{(\alpha_1,\alpha_2)}(A)$ with $\alpha_1>3/2$ and 
$\alpha_2>\max(\frac{\alpha_1}{2\alpha_1-3}, \frac{3\alpha_1}{2\alpha_1-1})$, then, under
assumptions of Theorem \ref{main}, for $n$ large enough, 
$$\E\|\pi\1_A-\tilde{\pi}\|_n^2\leq C\underset{m\in\mathcal{M}_n}{\inf}\left\{d^2(\pi\1_A,S_m) 
    +\frac{D_{m_1}D_{m_2}}{n}\right\}+\frac{C'}{n}.$$
\end{thm}

\begin{rem}
  The condition on the regularity of $\pi$ is verified for example if $\alpha_1>2$ and $\alpha_2>2$.
If $\alpha_1=\alpha_2=\alpha$, it is equivalent to $\alpha>2$.
\end{rem}

\emph{Proof: }
We recall that $\|\pi\|_\infty$ denotes actually $\|\pi\1_A\|_\infty$ and 
we introduce the following set:
$$\Lambda=\left\{\left|\frac{\|\hat\pi\|_\infty}{\|\pi\1_A\|_\infty}-1\right|<\frac{1}{2}\right\}.$$
As previously, we decompose the space:
$$\E\|\tilde{\pi}-\pi\1_{A}\|_n^2=
\E\left(\|\tilde{\pi}-\pi\1_{A}\|_n^2\1_{\Omega_{\rho}^*\cap\Lambda}\right)
+\E\left(\|\tilde{\pi}-\pi\1_{A}\|_n^2\1_{\Omega_{\rho}^*\cap\Lambda^c}\right)
+\E\left(\|\tilde{\pi}-\pi\1_{A}\|_n^2\1_{\Omega_{\rho}^{*c}}\right)$$
We have already dealt with the third term.
For the first term, we can proceed as in the proof of Theorem \ref{main} as soon as
$$\theta p(m,m')\leq \overline{\pen}(m)+\overline{\pen}(m')$$
with $\theta=3\rho=9/2$ and $p(m,m')=10\|\pi\|_\infty{D(m,m')}/{n}$.
But on $\Lambda$, $\|\pi\|_\infty<2\|\hat\pi\|_\infty$ and so
\begin{eqnarray*}
\theta p(m,m')&=&10\theta\|\pi\|_\infty\frac{D(m,m')}{n}
\leq 20\theta\|\hat\pi\|_\infty\frac{D(m,m')}{n}\\
&\leq& 20\theta\|\hat\pi\|_\infty\frac{D_{m_1}D_{m_2}}{n}
     +20\theta\|\hat\pi\|_\infty\frac{D_{m'_1}D_{m'_2}}{n}
\end{eqnarray*}
It is sufficient to set $\overline{K_0}=20\theta$.

Now, inequality \eqref{majdepi-pichapeau} gives
$$\E\left(\|\pi\1_A-\hat{\pi}_{\hat{m}}\|_n^2\1_{\Omega_{\rho}^{*}\cap\Lambda^c}\right)
\leq (\|\pi\|_\infty+4\phi_2 n^{1/3})P(\Omega_{\rho}^*\cap\Lambda^c).$$
It remains to prove that $P(\Omega_{\rho}^*\cap\Lambda^c)\leq Cn^{-4/3}$ for some constant $C$.
\begin{eqnarray*}\hspace{-0.8cm}
P(\Omega_{\rho}^*\cap\Lambda^c)&=&P(|\|\hat\pi\|_\infty-\|\pi\1_A\|_\infty|\1_{\Omega_{\rho}^{*}}
\geq\|\pi\|_\infty/2)\leq P(\|\hat\pi-\pi\1_A\|_\infty\1_{\Omega_{\rho}^{*}}\geq\|\pi\|_\infty/2)\\
&\leq&P(\|\hat\pi-\pi_{{m*}}\|_\infty\1_{\Omega_{\rho}^{*}}\geq\|\pi\|_\infty/4)
+P(\|\pi_{{m*}}-\pi\1_A\|_\infty\geq\|\pi\|_\infty/4)\\
&\leq&P\left(\|\hat\pi-\pi_{{m*}}\|\1_{\Omega_{\rho}^{*}}
\geq\frac{\|\pi\|_\infty}{4\phi_0\sqrt{D_{{m}_1*}D_{{m}_2*}}}\right)
+P(\|\pi_{{m*}}-\pi\1_A\|_\infty\geq\|\pi\|_\infty/4)\\
\end{eqnarray*}
since $\|\hat\pi-\pi_{{m*}}\|_\infty\leq\phi_0\sqrt{D_{{m*}_1}D_{{m*}_2}}
\|\hat\pi-\pi_{{m}^*}\|.$

Furthermore the inequality 
$\gamma_n(\hat{\pi})\leq \gamma_n(\pi_{m*})$ leads to
$$\|\hat{\pi}-\pi\1_A\|_n^2\leq \|\pi_{m*}-\pi\1_A\|_n^2
  +\frac{1}{\theta'}\|\hat{\pi}-\pi_{m*}\|_f^2+\theta'\sup_{t\in B_f({m*})}Z_n^2(t)$$
and then, on $\Omega_\rho$,
\begin{eqnarray*}\hspace{-0.4cm}
\|\hat{\pi}-\pi_{m*}\|_f^2\left(1-\frac{2\rho}{\theta'}\right)&\leq &4\rho\|\pi_{m*}-\pi\1_A\|_n^2
  +2\rho\theta'\sup_{t\in B_f({m*})}Z_n^2(t)\\
\text{ so }\qquad\|\hat{\pi}-\pi_{m*}\|^2&\leq &\frac{4\rho\theta'f_0^{-1}}{\theta'-2\rho}\|\pi_{m*}-\pi\1_A\|_n^2
  +\frac{2\rho\theta'^2f_0^{-1}}{\theta'-2\rho}\sup_{t\in B_f({m*})}Z_n^2(t)\\
&\leq &12\rho f_0^{-1}|A_2|\|\pi_{m*}-\pi\1_A\|_\infty^2
  +18\rho^2f_0^{-1}\sup_{t\in B_f({m*})}Z_n^2(t)\\
\end{eqnarray*}
with $\theta'=3\rho$ and by remarking that for $t$ with support $A$, 
$\|t\|_n^2\leq |A_2|\|t\|_\infty^2$.
Thus 
\begin{align}\nonumber
P(\Omega_{\rho}^*\cap\Lambda^c)&\leq P(\sup_{t\in B_f({m*})}Z_n^2(t)\1_{\Omega_{\rho}^{*}}
\geq\frac{\|\pi\|_\infty^2}{32\phi_0^2n^{1/3}}\frac{1}{18\rho^2f_0^{-1}})\\ \nonumber
&\quad+P(\|\pi_{{m*}}-\pi\1_A\|_\infty^2\geq\frac{\|\pi\|_\infty^2}{32\phi_0^2D_{m_1*}D_{m_2*}}
\frac{1}{12\rho f_0^{-1}|A_2|})\\ \nonumber
&\quad+P(\|\pi_{{m*}}-\pi\1_A\|_\infty\geq\|\pi\|_\infty/4)\\
\label{equ}
\begin{split}
&\leq P(\sup_{t\in B_f({m*})}Z_n^2(t)\1_{\Omega^{*}}\geq\frac{a}{n^{1/3}})
+P(D_{m_1*}D_{m_2*}\|\pi_{{m*}}-\pi\1_A\|_\infty^2\geq {b})\\ &\quad
+P(\|\pi_{{m*}}-\pi\1_A\|_\infty\geq\frac{\|\pi\|_\infty}{4})
\end{split}
\end{align}
with $a=\cfrac{\|\pi\|_\infty^2}{32\phi_0^2}\cfrac{1}{18\rho^2f_0^{-1}}$ and 
$b=\cfrac{\|\pi\|_\infty^2}{32\phi_0^2}\cfrac{1}{12\rho f_0^{-1}|A_2|}$.

We will first study the two last terms in \eqref{equ}.
Since the restriction $\pi_A$ of $\pi$ belongs to $B_{2,\infty}^{(\alpha_1,\alpha_2)}(A)$, the imbedding theorem proved in 
Nikol$'$ski{\u\i} \cite{Nikolskij} p.236 implies that $\pi_A$ belongs to $B_{\infty,\infty}^{(\beta_1,\beta_2)}(A)$
with $\beta_1=\alpha_1(1-1/\bar\alpha)$ and $\beta_2=\alpha_2(1-1/\bar\alpha)$.
Then the approximation lemma \ref{approx} (which is still valid for the trigonometric polynomial spaces 
with the infinite norm instead of the $L^2$ norm) yields to 
$$\|\pi_{m*}-\pi\1_A\|_\infty\leq C(D_{m_1*}^{-\beta_1}+D_{m_2*}^{-\beta_2}).$$
And then, since $D_{m_1*}=D_{m_2*}$,
\begin{eqnarray*}
 D_{m_1*}D_{m_2*}\|\pi_{m*}-\pi\1_A\|_\infty^2&\leq& C'(D_{m_1*}^{2-2\beta_1}+D_{m_1*}^{2-2\beta_2})\\
&\leq &C'((\ln n)^{2-2\beta_1}+(\ln n)^{2-2\beta_2})\rightarrow 0
\end{eqnarray*}
Indeed $\begin{cases}
2-2\beta_1<0 \Leftrightarrow 2\alpha_1\alpha_2-3\alpha_2-\alpha_1>0\\
2-2\beta_2<0 \Leftrightarrow 2\alpha_1\alpha_2-3\alpha_1-\alpha_2>0
\end{cases}$
and this double condition is ensured when  $\alpha_1>3/2$ and 
$\alpha_2>\max(\frac{\alpha_1}{2\alpha_1-3}, \frac{3\alpha_1}{2\alpha_1-1})$.
Consequently, for $n$ large enough,
$$P(D_{m_1*}D_{m_2*}\|\pi_{{m*}}-\pi\1_A\|_\infty^2\geq {b})
+P(\|\pi_{{m*}}-\pi\|_\infty\geq\frac{\|\pi\|_\infty}{4})=0.$$

We will now prove that 
$$P\left(\sup_{t\in B_f(m*)}Z_n^2(t)\1_{\Omega^{*}}\geq\frac{a}{n^{1/3}}\right)\leq \frac{C}{n^{4/3}}$$
and then using \eqref{equ}, we will have $P(\Omega_{\rho}^*\cap\Lambda^c)\leq {C}{n^{-4/3}}$.
We remark that, if $(\varphi_j\otimes\psi_k)_{j,k}$ is a base of $(S_{m*}, \|.\|_f)$,
$$\sup_{t\in B_f(m*)}Z_n^2(t)\leq \sum_{j,k}Z_n^2(\varphi_j\otimes\psi_k)$$
and we recall that, on $\Omega^*$, $Z_n(t)=\sum_{r=1}^4Z_{n,r}^*(t)$ (see the proof of
Proposition \ref{talagrand}). 
So we are interested in 
$$P\left(Z_{n,1}^{*2}(\varphi_j\otimes\psi_k)\1_{\Omega^{*}}\geq\frac{a}{4D_{m_1*}D_{m_2*}n^{1/3}}\right).$$

Let $x=B n^{-2/3}$ with $B$ such that $2f_0^{-2}B^2+4\|\pi\|_\infty B\leq a/4$ 
(for example $B=\inf(1,a/8(f_0^{-2}+2\|\pi\|_\infty)$). Then
$$(\sqrt{2\|\pi\|_\infty x}+\sqrt{D_{m_1*}D_{m_2*}}f_0^{-1}x)^2
\leq \frac{a}{4D_{m_1*}D_{m_2*}n^{1/3}}.$$

So we will now bound 
$P(Z_{n,1}^{*}(\varphi_j\otimes\psi_k)\1_{\Omega^{*}}\geq\sqrt{2\|\pi\|_\infty x}+
\sqrt{D_{m_1*}D_{m_2*}}f_0^{-1}x)$ 
by using the Bernstein inequality given in \cite{birge&massart98}.
That is why we bound $\E |\frac{1}{4q_n}\sum_{i=1, i \text{ odd}}^{2q_n-1} \Gamma_{i}^*(t)|^m$ 
for all integer $m\geq 2$,

\begin{eqnarray*}\hspace{-0.8cm}
\E |\frac{1}{4q_n}\sum_{i=1, i \text{ odd}}^{2q_n-1} \Gamma_{i}^*(t)|^m
&\leq& \cfrac{(2\|t\|_\infty q_n)^{m-2}}{(4q_n)^m}
\E\left|\sum_{i=1, i \text{ odd}}^{2q_n-1}[t(X_i^*,X_{i+1}^*)-\int t(X_i^*,y)\pi(X_i^*,y)dy]\right|^2\\
&\leq & \left(\cfrac{\|t\|_\infty}{2}\right)^{m-2}\cfrac{1}{16q_n^2}
\E\left|\sum_{i=1, i \text{ odd}}^{2q_n-1}[t(X_i,X_{i+1})-\int t(X_i,y)\pi(X_i,y)dy]\right|^2 \\
&\leq&\left(\cfrac{\|t\|_\infty}{2}\right)^{m-2}\cfrac{1}{16}
\int t^2(x,y)f(x)\pi(x,y)dxdy\\
&\leq& \frac{1}{2^{m+2}}(\|t\|_\infty)^{m-2}\|\pi\|_\infty\|t\|_f^2.
\end{eqnarray*}

Then 
\begin{eqnarray*}
\E |\frac{1}{4q_n}\sum_{i=1, i \text{ odd}}^{2q_n-1} \Gamma_{i*}(\varphi_j\otimes\psi_k)|^m
&\leq& \frac{1}{2^{m+2}}(\sqrt{D_{m_1*}D_{m_2*}}f_0^{-1})^{m-2}\|\pi\|_\infty.
\end{eqnarray*}
Thus the Bernstein inequality gives 
$$P(|Z_{n,1}^{*}(\varphi_j\otimes\psi_k)|\geq \sqrt{D_{m_1*}D_{m_2*}}f_0^{-1}x
+\sqrt{2\|\pi\|_\infty x})\leq 2e^{-p_nx}.$$
Hence \begin{eqnarray*}
P(\sup_{t\in B_f(m*)} Z_{n,1}^{*2}(t)\1_{\Omega^{*}}\geq \frac{a}{4n^{1/3}})
&\leq&2D_{m_1*}D_{m_2*}\exp\{-p_nB n^{-2/3}\}\\
&\leq&2n^{2/3}\exp\{-\frac{B}{4}\frac{n^{1/3}}{q_n} \}.
\end{eqnarray*}
But $2n^{2/3}\exp\{-\cfrac{B}{4}\cfrac{n^{1/3}}{q_n} \}\leq Cn^{-4/3}$ since $q_n\leq n^{1/6}$ and so 
$$P\left(\sup_{t\in B_f(m*)}Z_n^2(t)\1_{\Omega^{*}}\geq\frac{a}{n^{1/3}}\right)\leq \frac{4C}{n^{4/3}}.$$





\bibliographystyle{elsart-num}
\bibliography{biblio2}

\end{document}